\documentclass[a4paper,10pt]{amsart}
\usepackage{amsmath,amsthm,amssymb,amscd}
\usepackage{amsfonts}
\usepackage[american]{babel}
\usepackage[all]{xy}
\usepackage{fancyhdr}
\usepackage{enumerate}

\swapnumbers\newtheorem{prop}{Proposition}[section]
\newtheorem{theo}[prop]{Theorem}
\newtheorem{lem}[prop]{Lemma}
\newtheorem{coro}[prop]{Corollary}

 {\theoremstyle{definition} \newtheorem{de}[prop]{Definition} \theoremstyle{definition}
 \newtheorem{nota}[prop]{Notation} \newtheorem{des}[prop]{Definitions}
  \newtheorem{ex}[prop]{Example} \newtheorem{rem}[prop]{Remark}
 \newtheorem{exs}[prop]{Examples} \newtheorem{rems}[prop]{Remarks}
  }

{\swapnumbers\theoremstyle{definition}

}

\newcommand{\Q}{\mathbb{Q}}
\newcommand{\Z}{\mathbb{Z}}

\newcommand{\length}{\mathrm{length}}

\newcommand{\End}{\mathrm{End }}
\newcommand{\Hom}{\mathrm{Hom }}

\newcommand{\exseq}[3]{0\rightarrow {#1}\rightarrow{#2}\rightarrow{#3}\rightarrow 0}

\newcommand{\freealgebra}[2]{{#1} \langle {#2}\rangle}

\newcommand{\lexseqfivemap}[6]{0\rightarrow{#1}\rightarrow{#2} \stackrel{#6}{\rightarrow}{#3} \rightarrow \\\rightarrow{#4}\rightarrow{#5}}

\newcommand{\lexseqsixmap}[7]{0\rightarrow{#1}\rightarrow{#2} \stackrel{#7}{\rightarrow}{#3} \rightarrow \\\rightarrow{#4}\rightarrow{#5}\rightarrow{#6}}
\newcommand{\exseqmap}[4]{0\rightarrow {#1}\stackrel{#4}{\rightarrow}{#2}\rightarrow{#3}\rightarrow 0}
\newcommand{\lexseqmap}[4]{0\rightarrow {#1}\stackrel{#4}{\rightarrow}{#2}\rightarrow{#3}}

 \newenvironment{pro}{\noindent {\bf Proof.}\ }{\ \rule{1ex}{1ex}\vspace{0.5truecm}}

\DeclareMathOperator{\im}{im}
\DeclareMathOperator{\Ext}{Ext}\DeclareMathOperator{\Tor}{Tor}
\DeclareMathOperator{\lMod}{-Mod}\DeclareMathOperator{\rMod}{Mod-}
\DeclareMathOperator{\tr}{tr}\DeclareMathOperator{\pd}{pd}\DeclareMathOperator{\Gen}{Gen}
\DeclareMathOperator{\coker}{coker}
\DeclareMathOperator{\Add}{Add} 
\DeclareMathOperator{\leng}{length}

\DeclareMathOperator{\Tr}{Tr} \DeclareMathOperator{\Spec}{Spec}
\DeclareMathOperator{\ann}{ann}
\DeclareMathOperator{\maxspec}{max-spec}
\begin{document}
\title{Tilting modules arising from ring epimorphisms}

\author{Lidia Angeleri H\"ugel}
\address{Dipartimento di Informatica e
Comunicazione, Universit\`a degli Studi dell'Insubria, Via Mazzini
5, I - 21100 Varese, Italy} \email{lidia.angeleri@uninsubria.it}
\author{Javier S\'anchez }
\address{ Departament de Matem\`atiques \\
Universitat Aut\`onoma de Barcelona \\ E-08193 Be\-lla\-te\-rra
(Barcelona), Spain} \email{jsanchez@mat.uab.cat}


\date{\today}

\begin{abstract}
We show that a tilting module $T$ over a ring $R$ admits an exact sequence
 $\exseq R{T_0}{T_1}$ such that $T_0,T_1\in\Add(T)$ and
 ${\rm Hom}_R(T_1,T_0)=0$ if and only if $T$ has the form
 $S\oplus S/R$  for some injective ring epimorphism
 $\lambda:R\to S$ with the property that
$\Tor_1^R(S,S)=0$ and pd$S_R\le 1$. We then study the
case where $\lambda$ is a universal localization  in the
sense of Schofield \cite{Schofieldbook}. Using results from
  \cite{CrawleyBoeveytubes}, we give applications to tame
  hereditary algebras and hereditary noetherian prime rings.
  In particular, we show that every tilting module
  over a Dedekind domain or over a classical maximal order
   arises from universal localization.
\end{abstract}

\maketitle

\section*{Introduction}

Tilting theory started in the early eighties in representation theory
of finite dimensional algebras as a tool to relate two module categories
via functors inducing crosswise equivalences between certain parts
of both categories. Nowadays tilting plays an important role in
various branches of modern algebra, ranging from Lie theory and
algebraic geometry to homotopical algebra. We refer to \cite{Handbook} for a survey on such developments.

In this paper, we will  consider (large) tilting modules over an arbitrary ring $R$, according to the following definition.
A right $R$-module $T$ is said to be a \emph{tilting module} if it
satisfies the following properties:
\begin{enumerate}[(T1)]
\item $T$ has projective dimension at most one. \item
$\Ext_R^1(T,T^{(I)})=0$ for any set $I$. \item There exists an exact
sequence $\exseq R{T_0}{T_1}$ where $T_0,T_1$ are isomorphic to direct summands of direct sums of copies of $T$.
\end{enumerate}

A typical example of a (not finitely generated) tilting module is provided by the $\Z$-module $T=\Q\oplus \Q/\Z$.
 Its tilting class $\Gen T=\Gen\Q$ is the class  of all divisible groups.

In fact, following this pattern, one can use localization
techniques to construct tilting modules in many contexts. The
papers \cite{Angeleriherberatrlifaj2,Salce1,Salce2} already
contain  results in this direction. In the present paper, we push
forward this idea. We show that every injective  ring epimorphism
$\lambda:R\to S$ with the property that $\Tor_1^R(S,S)=0$ and $\pd
S_R\le 1$ gives rise to a tilting $R$-module $S\oplus S/R$
(Theorem \ref{theo:tiltingepimorphism}). Moreover, we characterize
the tilting modules that arise in this way. Namely, a tilting
module $T$ is equivalent to a tilting module $S\oplus S/R$ as
above if and only if the exact sequence
 $0\to R\to T_0\to T_1\to 0$ in condition (T3) can be chosen  with the additional
 property that ${\rm Hom}_R(T_1,T_0)=0$ (Theorem \ref{chara}).

Our construction yields many interesting examples of tilting
modules. For example, if the total ring of quotients
$Q^r_{\textrm{tot}}(R)$ of $R$ has projective dimension   one over
$R$, then $Q^r_{\textrm{tot}}(R)\oplus Q^r_{\textrm{tot}}(R)/R$ is
a tilting right $R$-module. Note that in general, however, the
ring epimorphism $\lambda:R\to S$ need not be a perfect
localization, see Examples \ref{notperfect}(1) and
 \ref{ex:freegroupalgebra}.  Examples of tilting modules that do not
 arise from  ring epimorphisms as above are given in
 Examples~\ref{notperfect}(3) and (4).

Given a tilting $R$-module $S\oplus S/R$ as above, in general it is difficult to compute its tilting class. In many cases, however,
 the tilting class can be described in terms of divisibility.

 For example, if $\mathfrak{U}$ is a left Ore set of
non-zero-divisors of $R$ such that  $\pd(\mathfrak{U}^{-1}R_R)\leq 1$,
then
$T_{\mathfrak{U}}=\mathfrak{U}^{-1}R\oplus \mathfrak{U}^{-1}R/R$
is a tilting right $R$-module whose tilting class $\Gen T_{\mathfrak{U}}$ coincides with the
 class of $\mathfrak{U}$-divisible right $R$-modules (Corollary \ref{coro:Oretiltingclass}).

 More generally,
if  $\mathcal{U}$ is a class of   finitely presented right
$R$-modules of projective dimension one such that $\Hom_R(\mathcal{U},R)=0$,
we can consider the universal localization $R_{\mathcal{U}}$ of $R$ at $\mathcal{U}$ in the sense of Schofield \cite{Schofieldbook}.
 Suppose that $R$ embeds in $R_{\mathcal{U}},$ and
$\pd (R_{\mathcal{U}})_R\leq 1$. Then
$T_{\mathcal{U}}=R_{\mathcal{U}}\oplus R_{\mathcal{U}}/R$ is a
tilting right $R$-module.
If we  further assume that
$R_{\mathcal{U}}/R$ is a direct limit of $\mathcal{U}$-filtered right $R$-modules,
 then the tilting class $\Gen T_{\mathcal{U}}$ coincides with the class $\mathcal{U}^\perp$ of all modules $M$ satisfying $\Ext^1_R(\mathcal{U},M)=0$
 (Corollary \ref{theo:tiltingwithclass}).

From  work of Schofield and Crawley-Boevey
\cite{CrawleyBoeveytubes} we know that universal
localizations satisfying such assumptions   occur  over
hereditary rings with a faithful rank function $\rho$
having the property that $R_\rho$ is simple artinian, see
Corollary \ref{coro:rhotorsiongivetilting}. An important
example for this situation is provided by
finite-dimensional tame hereditary algebras, see
\cite{CrawleyBoeveytubes}. In this case, we obtain  a
tilting module $T_\mathcal{U}=R_\mathcal{U}\oplus
R_\mathcal{U}/R$  with tilting class $\mathcal{U}^\perp$
 for every set
$\mathcal{U}$ of  simple regular
modules (Example \ref{ex:tame}). In a forthcoming paper, this result will be used to  classify the infinite-dimensional tilting modules over  tame hereditary algebras.

Another interesting example is the case of a hereditary noetherian prime ring with (simple artinian) quotient ring $A$.
In Theorem \ref{theo:tiltingprimenoetherianrings},
we prove that in this case  $T=A\oplus A/R$ is a tilting right $R$-module
with tilting class $\mathcal{U}^\perp$ where $\mathcal{U}$ is the class of all simple right $R$-modules.
Moreover, for any overring $R<S<A$ there exists a unique subset
$\mathcal{U}_S$ of $\mathcal{U}$ such that
$S\oplus S/R$ is a tilting right $R$-module with tilting
class $\mathcal{U}_S^\perp.$
And for any right Ore set
$\mathfrak{S}$ consisting of regular elements,
$T_{\mathfrak{S}}= R\mathfrak{S}^{-1}\oplus R\mathfrak{S}^{-1}/R$
is a tilting right $R$-module with tilting class
$\mathcal{U}_\mathfrak{S}^\perp$
where $\mathcal{U}_{\mathfrak{S}}$
is the class of all simple modules whose elements are annihilated by some element of
$\mathfrak{S}$.

In Corollary \ref{coro:alltiltingDedekind} we apply this result to
 a Dedekind domain $R$.
 We recover a classification result from \cite{BazzoniEklofTrlifaj}
 and show that the tilting modules $T_{\mathfrak{P}}=R_{\mathcal{U}_{\mathfrak{P}}}\oplus
R_{\mathcal{U}_{\mathfrak{P}}}/R$ arising from universal
localization at $\mathfrak{U}_\mathfrak{P}=\{R/\mathfrak{m}\mid
\mathfrak{m}\in\mathfrak{P}\}$, where $\mathfrak{P}$ runs through
all subsets of $\maxspec(R)$, form
 a representative set up to
equivalence of the class of all tilting $R$-modules.

\section{Homological properties of ring epimorphisms.}

\begin{nota}
For a ring $R$ (with 1), we denote by $\rMod R$ the category of all right $R$-modules.

Moreover, given $M\in\rMod R$, we write $\pd M$ for the projective
dimension of $M$.
\end{nota}

\begin{de}
Let $R,S$ be two rings. A morphism of rings $\lambda\colon
R\rightarrow S$ is called a \emph{ring epimorphism} if, for every
pair of  morphisms of rings $\delta_i\colon S\rightarrow T,\
i=1,2,$ the condition $\delta_1\lambda=\delta_2\lambda$ implies
$\delta_1=\delta_2.$
\end{de}

Of course, if  $\lambda\colon R\rightarrow S$ is a morphism of
rings, then every right (left) $S$-module is a right (left)
$R$-module, and every morphism of right (left) $S$-modules is a
morphism of right (left) $R$-modules. Moreover, it is well known
that   the category $\rMod S$ is a full subcategory of $\rMod R$
if and only if $\lambda$ is a ring epimorphism \cite[Chapter~XI,
Proposition~1.2]{Stenstrom}.

\bigskip

We will mainly deal with injective ring epimorphisms which
in addition satisfy the following homological property studied
 by  Geigle and Lenzing in \cite{GeigleLenzing}, see also
  \cite{DicksMayer-Vietoris} and \cite{NeemanNoncommutativeII}.

\begin{de}
Let $R,S$ be two rings and $\lambda\colon R\rightarrow S$ a ring
epimorphism. Then $\lambda$ is a  \emph{homological ring
epimorphism} if $\Tor_i^R(S,S)=0$ for all $i>0$.
\end{de}

Actually, we will see that in our context it is enough to require that $\Tor_1^R(S,S)=0$, a condition that Schofield has characterized as follows.

\begin{theo}[{\cite[Theorem~4.8]{Schofieldbook}}]\label{theo:epitor}
Let $\varphi\colon R\rightarrow S$ be a ring epimorphism. The
following statements are equivalent.
\begin{enumerate}[(1)]
 \item $\Tor_1^R(S,S)=0$.
 \item $\Tor_1^R(M,N)=\Tor_1^S(M,N)$ for all $M\in\rMod S$ and $N\in S\lMod$.
\item
$\Ext_R^1(M,N)=\Ext_S^1(M,N)$ for all $M,N\in \rMod S$.
\item $\Ext_R^1(M,N)=\Ext_S^1(M,N)$ for all $M,N\in S\lMod$.
\end{enumerate}
\end{theo}

\bigskip

The following notion from \cite{GeigleLenzing} will be useful for our discussion.

\begin{des}
If $\mathcal{S}$ is a class of right $R$-modules, the
\emph{(right) perpendicular category to $\mathcal{S}$} is defined
to be the full subcategory ${\mathcal X}_{\mathcal S}$ of $\rMod R$
consisting of all modules $A$ satisfying the following two
conditions
\begin{enumerate}[a)] \item $\Hom_R (S,A)=0$ for all
$S\in\mathcal{S}.$ \item $\Ext_R^1(S,A)=0$ for all
$S\in\mathcal{S}.$\end{enumerate}
\end{des}
If $\mathcal{S}=\{S\}$ we will write $\mathcal{X}_S$
instead of $\mathcal{X}_{\{S\}}.$

 Given an injective homological ring epimorphism $\lambda\colon R\rightarrow S$,  the right $S$-modules can be characterized inside the category of all right $R$-modules
as the objects of  the perpendicular category  ${\mathcal X}_{S/R}$.

\begin{theo}[{cf.
\cite[Proposition~4.12]{GeigleLenzing}}]\label{theo:Smodules} Let
$\lambda\colon R\rightarrow S$ be an injective ring epimorphism
with $\Tor_1^R(S,S)=0.$ Then the following are equivalent for
$M\in\rMod R$. \begin{enumerate}[(1)] \item $M\in\rMod S.$ \item
$\Ext_R^1(S/R,M)=\Hom_R(S/R,M)=0.$
\end{enumerate}
\end{theo}
\begin{pro}
Applying $\Hom_R({}_-,M)$ to the exact sequence $\exseq RS{S/R}$
we get
\begin{eqnarray*}
\lexseqfivemap{\Hom_R(S/R,M)}{\Hom_R(S,M)}{\Hom_R(R,M)}{\Ext_R^1(S/R,M)}{\Ext_R^1(S,M)}{\gamma}
\end{eqnarray*}
$(1)\Rightarrow (2):$ If $M\in\rMod S,$ then $\Ext_R^1(S,M)=\Ext_S^1(S,M)=0$.
Moreover the composition of maps
$M\cong\Hom_S(S,M)=\Hom_R(S,M)\stackrel{\gamma}{\rightarrow}\Hom_R(R,M)\cong
M$ is the identity on $M,$ and $\gamma$ is an isomorphism. Hence
$\Ext_R^1(S/R,M)=\Hom_R(S/R,M)=0.$
\\
$(2)\Rightarrow (1):$ Assume
$\Ext_R^1(S/R,M)=\Hom_R(S/R,M)=0$. Then $\gamma$ is an isomorphism, and $\Hom_R(S,M)\stackrel{\gamma}{\rightarrow}\Hom_R(R,M)\cong
M$, $f\mapsto f_{|R}\mapsto f(1)$, endows $M$ with a structure
of right $S$-module.
\end{pro}

\begin{rem}
 As a consequence of the last proof,  we see that for  a right $R$-module $M,$
the only possible structure   as right $S$-module is the one
given by $\Hom_R(S,M).$\\
\end{rem}


\section{Tilting modules arising from ring epimorphisms }\label{sec:tiltingfromepi}

\begin{des}\label{des:tiltingmodule} Let $R$ be a ring.

(1) Given a class $\mathcal{L}$ of right $R$-modules, we denote $$\mathcal{L}^\perp=\{M\in\rMod R\mid
\Ext_R^1(L,M)=0 \textrm{ for all } L\in\mathcal{L}\}.$$ If
$\mathcal{L}=\{L\}$ we will write $L^\perp$ instead of
$\{L\}^\perp.$

(2) For a right $R$-module $M,$ we denote by $\Add M$   the class of all isomorphic images of direct summands of
direct sums of copies of $M$, and by $\Gen M$  the class of all
right $R$-modules generated by $M,$ i.e. the right $R$-modules which
are epimorphic images of arbitrary direct sums of copies  of $M.$

(3)
A right $R$-module $T$ is said to be a \emph{tilting module} if  $\Gen
T=T^\perp.$  This is equivalent to the definition given in the introduction, see
  \cite{ColpiTrlifaj}. The class $T^\perp$ is called a \emph{tilting class}.

(5) Two tilting modules $T$ and $T'$ are said to be \emph{equivalent} if their tilting classes $T^\perp$ and $T'\,^\perp$ coincide.
This is equivalent to the condition $\Add(T)=\Add(T')$.
\end{des}

\begin{ex}
$\mathbb{Q}\oplus \mathbb{Q}/\mathbb{Z}$ is a tilting
$\mathbb{Z}$-module. Its tilting class is the class
$\Gen\Q$ of all divisible groups. Notice that
$\mathbb{Z}\hookrightarrow \mathbb{Q}$ is a ring
epimorphism with $\pd \mathbb{Q_Z}\leq 1$ and
$\Tor_1^\mathbb{Z}(\mathbb{Q},\mathbb{Q})=0.$
\end{ex}

We will now  study tilting
modules, like $\mathbb{Q}\oplus \mathbb{Q}/\mathbb{Z}$,  constructed
 from injective ring epimorphisms.
We start out with  a generalization of some results from
\cite[Section~6]{Angeleriherberatrlifaj2}, which in turn generalized part of
 \cite[Chapter~1]{Matliscotorsion}.
The proofs are very
similar to the ones in \cite{Angeleriherberatrlifaj2}, so we mostly omit them.

\begin{lem}\label{lem:traceepi}
Let $\lambda\colon R\rightarrow S$ be an injective ring
epimorphism, and let $M$ be a right $R$-module. The image of the
morphism $\Hom_R(S,M)\rightarrow M,\, f\mapsto f(1)$ coincides
with the trace $\tr_S(M)=\sum \{ f(S)\mid f\in\Hom_R(S,M)\}$ of
$S$  in $M$.
\end{lem}

\medskip

The following lemma  generalizes
 \cite[Lemma~1.8]{Matliscotorsion}.

\begin{lem}\label{lem:traceextension}
Let $R$ be a  ring. Let $\lambda\colon R\rightarrow S$ be an
injective ring epimorphism with $\Tor_1^R(S,S)=0$. Then the
following statements are equivalent.
\begin{enumerate}[a)]
\item $\tr_S(M/\tr_S(M))=0$ for all $M\in\rMod R$.
\item $\Gen S_R$ is closed under extensions.
\end{enumerate}
\end{lem}
\begin{pro}
$a)\Rightarrow b):$ Let $\exseq AB{B/A}$ be an exact sequence with
$A,B/A\in\Gen S_R.$ Since $A$ is contained in $\tr_S(B),$ we get
the surjective morphism of right $R$-modules $B/A\rightarrow
B/\tr_S(B).$ Hence $B/\tr_S(B)\in\Gen S_R,$ but by hypothesis
$\Hom_R(S,B/\tr_S(B))=0.$ Therefore $B/\tr_S(B)=0$ and
$B=\tr_S(B)\in\Gen S_R.$

$b)\Rightarrow a):$ Suppose $\tr_S(M/\tr_S(M))\neq0$ for some
right $R$-module $M.$ Then there exists a submodule $X$ of $M$
such that $X$ contains $\tr_S(M),$ $X/\tr_S(M)\neq 0$ and
$X/\tr_S(M)\in\Gen S_R.$ Consider the exact sequence $\exseq
{\tr_S(M)} X {X/\tr_S(M)}.$ By hypothesis, $X\in \Gen S_R,$ which
implies $X=\tr_S(M),$ a contradiction.
\end{pro}

\begin{theo}\label{theo:tiltingepimorphism}
Let $R$ be a ring. Let $\lambda\colon R\rightarrow S$ be an
injective ring epimorphism with $\Tor_1^R(S,S)=0$. Denote by
${\mathcal X}_S$ the perpendicular category  to $S_R$. The
following conditions are equivalent.
\begin{enumerate}[(1)]
\item $\pd(S_R)\leq 1$. \item  ${\mathcal X}_S$   is closed under cokernels of monomorphisms. \item
$\Ext_R^1(S/R, M)$ belongs to  ${\mathcal X}_S$   for any right $R$-module
$M$. \item $(S/R)^\perp=\Gen S_R$. \item $T=S\oplus S/R$ is a
tilting right $R$-module. \item $\pd((S/R)_R)\leq 1$.
\end{enumerate}
Moreover, under (1)-(6),
$\Hom_R(S,M/\tr_S(M))=0$ for any right $R$-module $M$.
\end{theo}

\begin{pro}
For $(1)\Rightarrow (2)$, see for example \cite[Proposition
1.1]{GeigleLenzing} or
\cite[Proposition~6.3]{Angeleriherberatrlifaj2}.

\noindent $(2)\Rightarrow (3):$ First of all, using that
$\lambda\colon R\rightarrow S$ is an injective ring
epimorphism with $\Tor_1^R(S,S)=0$, one shows as in
\cite[Lemma~6.1]{Angeleriherberatrlifaj2} that
$\Hom_R(S/R,M)$ belongs to ${\mathcal X}_S$ for any right
$R$-module $M.$ Denote now by $E(M)$ the injective hull of
$M.$ Applying the functor $\Hom_R(S/R,{}_-)$ to the exact
 sequence $\exseq M{E(M)}{E(M)/M}$ we get \begin{eqnarray*}
\lexseqfivemap
{\Hom_R(S/R,M)}{\Hom_R(S/R,E(M))}{\Hom_R(S/R,E(M)/M)}{\Ext_R^1(S/R,M)}{\Ext_R^1(S/R,E(M))=0}{\beta}.
 \end{eqnarray*}
Our assumption $(2)$ then yields that $\im\beta$ and
$\Ext_R^1(S/R,M)$ belong to ${\mathcal X}_S$.

\noindent $(3)\Rightarrow (4):$ Let $M$ be a right $R$-module.
Applying $\Hom_R({}_-,M_R)$ to the sequence $\exseq RS{S/R}$ we
obtain \begin{eqnarray*}
\lexseqsixmap{\Hom_R(S/R,M)}{\Hom_R(S,M)}{\Hom_R(R,M)}{\Ext_R^1(S/R,M)}{\Ext_R^1(S,M)}{\Ext_R^1(R,M)=0}{\alpha}.
\end{eqnarray*}
The natural isomorphism $\Hom_R(R,M)\rightarrow M$, defined by
$f\mapsto f(1)$, gives a map $\alpha\colon \Hom_R(S,M)\rightarrow
M$ whose image is the trace of $S$ in $M$ by
Lemma~\ref{lem:traceepi}. Hence $M\in\Gen S_R$ if and only if
$\alpha$ is surjective. If $M\in(S/R)^\perp,$ then clearly $\alpha$ is
surjective and $M\in\Gen S_R.$ Conversely, suppose that $\alpha$ is
surjective. Then $\Ext_R^1(S/R,M)\cong\Ext_R^1(S,M)$, so
$\Ext_R^1(S,M)$ belongs to ${\mathcal X}_S$ by $(3)$. But $\Ext_R^1(S,M)$ is a
right $S$-module, and the only right $S$-module which
belongs to ${\mathcal X}_S$ is the zero module. Hence
$\Ext_R^1(S,M)=\Ext_R^1(S/R,M)=0$, and $M\in(S/R)^\perp.$

\noindent $(4)\Rightarrow (5):$  By $(4),$
$\Gen T_R=\Gen S_R=(S/R)^\perp=T^\perp,$ and so $T$ is a tilting module.

\noindent $(5)\Rightarrow (6):$ If $T_R$ is a tilting right
$R$-module, then $\pd T_R\leq1$, which clearly implies $\pd
((S/R)_R)\leq 1$.

\noindent $(6)\Rightarrow (1)$ is clear.

To prove the last part of the Theorem, notice that  $\Gen
S_R$ is closed under extensions by $(4)$. Now apply
Lemma~\ref{lem:traceextension}.
\end{pro}

\begin{rems}\label{rems:isageneralizationofAHT}
Suppose $\lambda\colon R\rightarrow S$ is a morphism of rings as
in Theorem~\ref{theo:tiltingepimorphism}.
\begin{enumerate}[(1)]
\item
 When $R$ is a commutative ring and $S$ the full ring of quotients of $R,$
 the  objects of the right
perpendicular category ${\mathcal X}_S$  are precisely the
$R$-modules that Matlis called \emph{cotorsion} in
\cite{Matliscotorsion}.

\item  In many cases, for
example if $R$ is a hereditary ring, $S\oplus S/R$ is a two-sided
tilting $R$-module. \end{enumerate}
\end{rems}

\begin{exs}\label{coro:ringofquotients}
Let $R$ be a ring.
\begin{enumerate}[(1)]
\item Denote by $Q_{\textrm{max}}^r(R)_R$ the \emph{maximal right
ring of quotients} of $R$, see e.g. \cite[p. 200]{Stenstrom}.
Assume that $\pd(Q_{\textrm{max}}^r(R)_R)\leq 1$ and that one of
the following conditions is satisfied:
\begin{enumerate}[(a)]
\item $R$ is a right nonsingular ring such that every finitely
generated non-singular right $R$-module can be embedded in a free
module,
or
\item $Q_{\textrm{max}}^r(R)$ is right Kasch (for
example, this holds true whenever $Q_{\textrm{max}}^r(R)$ is semisimple).
\end{enumerate}
Then $Q_{\textrm{max}}^r(R)\oplus Q_{\textrm{max}}^r(R)/R$ is a
tilting right $R$-module.
This follows combining Theorem~\ref{theo:tiltingepimorphism} with
\cite[Chapter~XII,~Theorem~7.1]{Stenstrom} in case (a), or with
\cite[Chapter~XI,~Proposition~5.3]{Stenstrom} in case (b).

\smallskip

\item
By
\cite[Chapter~XI,~Theorem~4.1]{Stenstrom}
there exist  a ring
$Q^r_{\textrm{tot}}(R)$ and a ring epimorphism $\varphi\colon
R\rightarrow Q^r_{\textrm{tot}}(R)$ such that
\begin{enumerate}[i)] \item $\varphi$ is an injective ring epimorphism
and $Q^r_{\textrm{tot}}(R)$ is flat as a left $R$-module. \item
For every injective epimorphism $\gamma\colon R\rightarrow T$ of
rings such that ${}_RT$ is flat, there is a unique morphism of
rings $\delta\colon T\rightarrow Q^r_{\textrm{tot}}(R)$ such that
$\delta\gamma=\varphi$.
\end{enumerate}
 If $\pd (Q^r_{\textrm{tot}}(R)_R)\leq 1,$ then we infer from
Theorem~\ref{theo:tiltingepimorphism} that
$Q^r_{\textrm{tot}}(R)\oplus Q^r_{\textrm{tot}}(R)/R$ is a tilting
right $R$-module.
\end{enumerate}
\end{exs}

\medskip

Our next aim is to characterize the tilting modules that arise
from injective ring epimorphisms as in Theorem~\ref{theo:tiltingepimorphism}.
 We first introduce some terminology.

\begin{des}Let $R$ be a ring.
Let $M$ be a right $R$-module and $\mathcal{C}$ a class of right
$R$-modules closed under isomorphic images.

(1)
$\mathcal{C}$ is said to be a \emph{torsion class} if it is
closed under extensions, direct sums and epimorphic images.

(2)
  A morphism
$f\in\Hom_R(M,C)$ with $C\in\mathcal{C}$ is a
$\mathcal{C}$-\emph{preenvelope} of $M$ provided the morphism of
abelian groups
$\Hom_R(f,C')\colon\Hom_R(C,C')\rightarrow\Hom_R(M,C')$  is
surjective for each $C'\in\mathcal{C},$ that is, for each morphism
$f'\colon M\rightarrow C'$ there is a morphism $g\colon
C\rightarrow C'$ such that the following diagram is commutative.
\begin{eqnarray*}
\xymatrix@C=0.8cm@R=0.8cm{ M\ar[r]^f\ar[dr]_{f'} & C\ar@{-->}[d]^g
\\ & C'}
\end{eqnarray*}

(3)
A $\mathcal{C}$-preenvelope $f\in\Hom_R(M,C)$ is a
$\mathcal{C}$-\emph{envelope} of $M$ provided that $f$ is
\emph{left minimal}, that is,   every $g\in\End_R(C)$ such that
$f=gf$  is an automorphism.

(4) A $\mathcal{C}$-preenvelope $f\in\Hom_R(M,C)$ is said to be a
$\mathcal{C}$-\emph{reflection} of $M$  provided the morphism of
abelian groups
$\Hom_R(f,C')\colon\Hom_R(C,C')\rightarrow\Hom_R(M,C')$ is
bijective for each $C'\in\mathcal{C},$ that is, the morphism
$g\colon C\rightarrow C'$ in the diagram above is always uniquely
determined. Of course, every $\mathcal{C}$-{reflection} is a
$\mathcal{C}$-{envelope}.

(5)
 $\mathcal{C}$ is said to be a \emph{reflective} subcategory of $\rMod R$ if every $R$-module $M$ admits a $\mathcal{C}$-reflection.
\end{des}

\begin{rems}\label{refle}
(1) If $T$ is a tilting module, then every $R$-module $M$ admits a
$T^\perp$-preenvelope, see \cite{ATT}.  In particular, if  $0\to
R\stackrel{a}{\to} T_0\to T_1\to 0$ is an exact sequence where
$T_0,T_1\in \Add T$ as in condition (T3) in the Introduction, then
the map $a:R\to T_0$ is a $T^\perp$-preenvelope of $R$.

\smallskip

(2) It is well known that a class of right $R$-modules
$\mathcal{C}$ is a reflective subcategory of $\rMod R$ if and only
if the inclusion functor $\iota: \mathcal{C}\to \rMod R$ has a
left adjoint $\ell:\rMod R\to \mathcal{C}$. Then a $\mathcal
C$-reflection of $M$ is given as $\epsilon_M:M\to \ell(M)$ where
$\epsilon:1_{\rMod R}\to \iota\,\ell$ is the unit of the
adjunction, see e.g. \cite[Chapter X, \S 1]{Stenstrom}.

\smallskip

(3) When they exist, $\mathcal{C}$-envelopes and $\mathcal{C}$-reflections are uniquely determined up to isomorphism.
\end{rems}

\bigskip

We are now ready for the main result of this section.

\begin{theo}\label{chara} Let $R$ be a   ring and $T$ be a tilting right $R$-module. The following statements are equivalent.
\begin{enumerate}[(1)]
\item There is an injective ring epimorphism $\lambda:R\to S$ such
that ${\rm Tor}^R_1(S,S)=0$ and $S\oplus S/R$ is a tilting module
equivalent to $T$. \item There is an exact sequence $0\to
R\stackrel{a}{\to} T_0\to T_1\to 0$ such that $T_0,T_1\in \Add T$
and ${\rm Hom}_R(T_1,T_0)=0$.
\end{enumerate}
Moreover, under these conditions, $a:R\to T_0$ is a $T^\perp$-envelope of $R$, and $\lambda:R\to S$ is a homological ring epimorphism.
\end{theo}
\begin{pro}
The implication (1)$\Rightarrow$(2) follows immediately by choosing the exact sequence $0\to R\to S\to S/R\to 0$, keeping in mind that $\Add(S\oplus S/R)=\Add(T)$, and that  ${\rm Hom}_R(S/R,S)=0$ by Theorem \ref{theo:Smodules}.

For the implication (2)$\Rightarrow$(1), observe first that $a$ is
a $\Gen(T)$-preenvelope of $R$, so $\Gen T= \Gen T_0$ by
\cite[Lemma~1.1]{ATT}. Moreover, it is easy to see that
$T'=T_0\oplus T_1\in \Add T$   is a tilting module with $\Gen
T'=\Gen T_0$. Then $T'$ is  equivalent to $T$, and $\Gen
T=(T_0\oplus T_1)^\perp =(T_1)^\perp$. In particular, $T_1^\perp$
is a torsion class, and $T_1$ is a partial tilting module in the
sense of \cite{CTT}.

Denote now by $\mathcal{X}= \mathcal{X}_{T_1}$ the perpendicular
category of $T_1$.  As shown in \cite[Proposition 1.3 and
1.4]{CTT}, $\mathcal{X}$ is a reflective   subcategory of Mod$R$
which is closed under extensions, arbitrary direct sums and direct
products,  kernels and cokernels.

Then, as in \cite[Proposition 1.5]{CTT}, one can apply results of
Gabriel and de la Pe\~na \cite[Theorem~1.2]{GP} or Geigle and
Lenzing \cite[Proposition~3.8]{GeigleLenzing} to obtain a ring
epimorphism $\lambda:R\to S$ such that  the category $\rMod S$,
when viewed as a full subcategory of $\rMod R$,  is equivalent to
$\mathcal{X}$. More precisely, if  $\ell:\rMod R\to \mathcal{X}$
is a left adjoint of the inclusion functor $\iota: \mathcal{X}\to
\rMod R$, then $\ell(R)$ is a projective generator of
$\mathcal{X}$, and the functor $\Hom_R(\ell(R),-)$ preserves
coproducts. So, if we set $S=\End_R\ell(R)$,  we obtain mutually
inverse  functors $\Hom_R(\ell(R),-)$ and $-\otimes_S\ell(R)$
between  $\mathcal{X}$ and  $\rMod S$. Moreover, the assignment
  $\lambda (r)=\ell(m_r)$, where $m_r:R\to R$ denotes the left multiplication with the  element $r$, defines a ring epimorphism $\lambda: R\to S$. Note that
$\lambda (r)$ is
the uniquely determined element of $S$ extending   the endomorphism $m_r$ of $R$. In other words, if $\epsilon_R:R\to \ell(R)$ is the $\mathcal{X}$-reflection of $R$, then  $\lambda(r)\, \epsilon_R=\epsilon_R\, m_r$.

 Now,  in our case $T_0$ belongs to ${\mathcal X}$, and $R\stackrel{a}{\to} T_0$
 is an
  ${\mathcal X}$-reflection
  of $R$, that is, we can choose $\epsilon_R=a$ and $\ell(R)=T_0$.
   In particular, since  $a:R\to T_0$ is a monomorphism,
  $\lambda:R\to S$   is injective. Indeed,  $\lambda(r)=0$ implies $a\, m_r=0$, hence $m_r=0$ and $r=0$.
Moreover, $T_0$ and $S_R$ are canonically isomorphic, thus $S_R$
has projective dimension at most one, and $\lambda$ is  a
homological ring epimorphism by
\cite[Corollary~4.8]{GeigleLenzing}. So, we conclude from Theorem
\ref{theo:tiltingepimorphism} that $S\oplus S/R$ is a tilting
module with tilting class $\Gen(S_R)=\Gen(T_0)=\Gen(T)$.

For the last statement, note that the $\mathcal{X}$-reflection $a:R\to T_0$ is  left minimal, thus also a  Gen$T$-envelope of $R$. \end{pro}

\begin{exs}\label{notperfect}
(1) In general, in the situation of Theorem \ref{chara}, $S$ is not flat as a left $R$-module.

Let us look at \cite[Example 2.2]{CTT}. Here $R$ is   the path algebra given by the quiver $3\to 1\leftarrow 2\leftarrow 4$, and we consider
  the tilting $R$-module $T=\tau I_1\oplus I_2\oplus I_3\oplus I_4$.
  There is an exact sequence $0\to R\to (\tau I_1)^3\oplus (I_4)^3\to I_2^4\oplus I_3\to 0$ where $T_1=I_2^4\oplus I_3$ and
  $T_0=(\tau I_1)^3\oplus (I_4)^3$ belong to $\Add(T)$, and
    ${\rm Hom}_R(T_1,T_0)=0$. But the left adjoint
     $\ell:\rMod R\to \mathcal{X}$ of the inclusion functor
      $\iota: \mathcal{X}\to \rMod R$ is not (left) exact,
      and thus $_RS$ is not perfect, cf.~\cite[Example~2.2]{CTT}.

An example where $S$ is not flat as a right (nor as a left)
$R$-module will be given in Example~\ref{ex:freegroupalgebra}.

\smallskip

(2)
If we omit the assumption ${\rm Hom}_R(T_1,T_0)=0$,
we still have a ring epimorphism $\lambda:R\to S$ such
 that $\mathcal{X}$ is equivalent to
$\rMod S$. However, $S_R$ need not have projective
dimension at most one, and $\lambda$ need not be a homological
 epimorphism. Consider \cite[Example 2.4]{CTT}. Here $R$ is
 the algebra given by the quiver
 $1\stackrel{\alpha}{\leftarrow} 2\stackrel{\beta}{\leftarrow} 3$
 with the relation $\beta\alpha=0$, and we take $T=R$ with the
 (split) exact sequence $0\to R\to R\oplus P_2\to P_2\to 0$.
  Then   $S_R\cong S_3\oplus S_1$ has projective dimension 2 and
  $\lambda$ is not a homological epimorphism,
  cf.~\cite[p.229]{CTT}.

\smallskip

(3) Let $R$ be a  hereditary (indecomposable) artin algebra of
infinite representation type. Denote by $\bf p$ the preprojective
component of $R$. There is a countably infinitely generated
tilting  $R$-module generating ${\bf p}^\perp$, called the {\em
Lukas tilting module}, and denoted by $L$, cf.\ \cite{L1,KT}. It
has the property that there are non-zero morphisms between any two
non-zero modules from $\Add{L}$, see \cite[Theorem~6.1 (b)]{L1}
and \cite[Lemma~3.3 (a)]{Lukasinfinite-rankmodules}.
 So, there cannot be
an exact sequence $0\to R\stackrel{a}{\to} L_0\to L_1\to 0$ such
that $L_0,L_1\in \Add L$ and ${\rm Hom}_R(L_1,L_0)=0$, and
therefore $L$ does not arise from a ring epimorphism as above.

\smallskip

(4) Let $R$ be a Pr\"ufer domain which is not a Matlis domain,
that is, the quotient field $Q$ has projective dimension $>1$ over
$R$. Then $R$ has no divisible envelope, see
\cite[Corollary~6.3.18]{GobelTrlifaj}. So, the Fuchs tilting
module $\delta$, which is a tilting module generating the class of
all divisible modules \cite[Example~5.1.2]{GobelTrlifaj}, is
another example of a tilting module  that does not arise from a
ring epimorphism as above.
\end{exs}

\section{Tilting modules arising from universal localization}

Let us recall Schofield's notion of universal localization.

\begin{nota}
Let $R$ be a ring. By $\mathcal{P}_R$ (${}_R\mathcal{P}$) we
denote the category of all finitely generated projective right
(left) $R$-modules.

Let $P$ and $Q$ be finitely generated projective right
$R$-modules. By $P^*$ we denote the finitely generated projective
left $R$-module $\Hom_R(P,R).$ If $\alpha\in\Hom_R(P,Q),$ we
denote by $\alpha^*$ the morphism of finitely generated left
$R$-modules $\alpha^*\colon Q^*\rightarrow P^*$ defined by
$\gamma\mapsto\gamma\alpha.$

For a class of morphisms  $\Sigma$  between finitely generated
projective right $R$-modules we set $\Sigma^*=\{\alpha^*\mid \alpha\in \Sigma\}.$

\end{nota}

\begin{theo}[{\cite[Theorem~4.1]{Schofieldbook}}]\label{def:universallocalization}
Let $R$ be a ring and $\Sigma$ be a class of morphisms between
finitely generated projective right $R$-modules. Then
there are a ring $R_\Sigma$ and a morphism of rings
$\lambda\colon R\rightarrow R_\Sigma$ such that
\begin{enumerate}[(i)]
\item $\lambda$ is \emph{$\Sigma$-inverting,} i.e. if
$\alpha\colon P\rightarrow Q$ belongs to  $\Sigma$, then
$\alpha\otimes_R 1_{R_\Sigma}\colon P\otimes_R R_\Sigma\rightarrow
Q\otimes_R R_\Sigma$ is an isomorphism of right
$R_\Sigma$-modules, and \item $\lambda$ is \emph{universal
$\Sigma$-inverting}, i.e. if $S$ is a ring such that there exists
a $\Sigma$-inverting morphism $\psi\colon R\rightarrow S$, then
there exists a unique morphism of rings $\bar{\psi}\colon
R_\Sigma\rightarrow S$ such that $\bar{\psi}\lambda=\psi$.
\end{enumerate}
\end{theo}

\begin{de} $\lambda\colon R\rightarrow R_\Sigma$
as above is called the \emph{universal localization of $R$ at
$\Sigma$}. In the same way one defines the universal localization
at a class of morphisms between finitely generated projective left
$R$-modules.
\end{de}

\begin{rems}\label{rems:universalisepi} Let $R$ be a ring and let $\Sigma$ be a class of
morphisms between finitely generated projective right $R$-modules.
\begin{enumerate}[(1)]
\item The universal localization $R_\Sigma$ of $R$ at $\Sigma$ is
unique up to isomorphism of $R$-rings, i.e. if $\lambda_i\colon
R\rightarrow S_i,\ i=1,2,$ are universal localizations of $R$ at
$\Sigma,$ there exists a unique isomorphism of rings
$\varphi\colon S_1\rightarrow S_2$ such that
$\varphi\lambda_1=\lambda_2.$

\item
$R_{\Sigma^*}$ is isomorphic to $R_\Sigma.$

\item $\lambda\colon R\rightarrow R_\Sigma$ is a ring epimorphism
with  $\Tor_1^R(R_\Sigma,R_\Sigma)=0.$ So, if $\lambda$ is
injective and $\pd R_\Sigma\le 1$, then we infer from Theorem
\ref{theo:tiltingepimorphism} that $ R_\Sigma\oplus  R_\Sigma/R$
is a tilting module with tilting class $\Gen R_\Sigma=
(R_\Sigma/R)^\perp$.

\item
If $\lambda\colon
R\rightarrow R_\Sigma$ is injective, then every $\alpha\in\Sigma$
is injective.

\end{enumerate}
\end{rems}
\begin{pro}
\noindent  (1) follows from the universal property of universal localization. For a proof of (2), we refer to
 \cite[pages 51-52]{Schofieldbook}. (3) holds   by \cite[Theorem~4.7]{Schofieldbook} and
Theorem~\ref{theo:epitor}.   For $(4)$ see  \cite[Proposition~2.2]{NeemanNoncommutativeII}.
\end{pro}

 Remarks \ref{rems:universalisepi} (2) and (4) show that the injectivity of all morphisms in $\Sigma$ and $\Sigma^*$ is a necessary condition for $\lambda$ to be injective. Let us now turn to the cokernels of maps in $\Sigma$.

\begin{des}
(a) Let $U$ be a right (left) $R$-module. We say that $U$ is a
\emph{bound right (left) $R$-module} if $U$ is finitely presented,
$\pd U=1$ and $\Hom_R(U,R)=0.$ In other words, $U$ is a bound
right (left)  $R$-module if and only if $U$ is the cokernel of
some morphism $\alpha\colon P\rightarrow Q$ with
$P,Q\in\mathcal{P}_R$ (${}_R\mathcal{P}$) such that $\alpha$ and
$\alpha^*$ are injective.

\smallskip

(b) If $U$ is  a bound module with projective presentation $\exseqmap PQU\alpha$,
  then we have an
exact sequence $\exseqmap{Q^*}{P^*}{\coker\alpha^*}{\alpha^*}$,
{and} $\coker\alpha^*$ is the \emph{Auslander-Bridger transpose}
of $U$ denoted by $\Tr U=\coker\alpha^*$,  see for example
\cite{AuslanderReitenSmalo}.

\smallskip

(c) Let now $\mathcal{U}$ be a class of  bound right $R$-modules.
For each $U\in\mathcal{U},$ consider a morphism $\alpha_U$ between
finitely generated projective right $R$-modules such that
\begin{equation}\label{eq:bounddefinition}
\exseqmap PQU{\alpha_U}\end{equation}
 is exact. We will denote by
$R_{\mathcal{U}}$ the universal localization of $R$ at
$\Sigma=\{\alpha_U\mid U\in\mathcal{U}\}.$ In fact,   $R_{\mathcal
U}$ does not depend on the chosen class $\Sigma$,
cf.~\cite[Theorem~0.6.2]{Cohnfreerings}, and we will also call it
the \emph{universal localization of $R$ at ${\mathcal{U}}$}. By
abuse of notation,  we will
  write $\alpha_U\in\mathcal{U}$ for any morphism
$\alpha_U$ between finitely generated projective right $R$-modules
as in \eqref{eq:bounddefinition} with $U\in\mathcal{U}$.

 Finally, a  right
$R$-module $N$  is said to be \emph{$\mathcal{U}$-torsion-free} if
$\Hom_R(U,N)=0$ for all $U\in\mathcal{U},$ and $N$  is said to be
\emph{$\mathcal{U}$-divisible} if $\Ext_R^1(U,N)=0$ for all
$U\in\mathcal{U}.$
\end{des}

\begin{rem}\label{rems:boundmodules}
It is known that a tilting module $T$ is of \emph{finite
type} \cite[Theorem~2.6]{BazzoniHerbera}, that is, there
exists a set $\mathcal{V}$ of finitely presented modules of
projective dimension at most one such that
$T^\perp=\mathcal{V}^\perp$. When $R$ is a semihereditary
ring, every finitely presented module $M$ is of the form
$M=P\oplus U$ where $P$ is finitely generated projective
and $U$ is a bound module
\cite[Theorem~1.2(3)]{LuckHilbertmodules}. Thus every
tilting class is of the form $\mathcal{V}^\perp$ where
$\mathcal{V}$ is a set of bound modules if $R$ is a
semihereditary ring.
\end{rem}

\begin{ex}\label{ex:oreisuniversal}
If $R$ is a ring and $\mathfrak{U}\subset R$ is a right
denominator set, then the right  Ore localization
$R\mathfrak{U}^{-1}$  is the universal localization of $R$ at all
the maps $\alpha_u\colon R\rightarrow R,\,r\mapsto
ur,$ where $u\in \mathfrak{U}.$ Equivalently,
$R\mathfrak{U}^{-1}$ is the universal localization of $R$ at the
maps $\alpha_u^*$,  $u\in\mathfrak{U}$, given by right
multiplication by $u.$

In the same way, if $\mathfrak{U}$ is a left denominator set, then
$\mathfrak{U}^{-1}R$ is the universal localization of $R$ at the
maps $\alpha_u,$ $u\in\mathfrak{U}$, and
also the universal localization at the maps
$\alpha_u^*$, $u\in\mathfrak{U}.$

Notice further that if $\mathfrak{U}$ is a left Ore set of
non-zero-divisors of $R,$ then $\mathcal{U}=\{R/uR\mid
u\in\mathfrak{U}\}$ and $\Tr\mathcal{U}=\{\Tr U \mid U\in\mathcal{U}\}=\{R/Ru\mid
u\in\mathfrak{U}\}$ are sets of
bound modules, and
$\mathfrak{U}^{-1}R=R_\mathcal{U}=R_{\Tr\mathcal{U}}.$ Moreover,
$M$ is $\mathfrak{U}$-torsion-free  iff $M$ is
$\mathcal{U}$-torsion-free, and $M$ is $\mathfrak{U}$-divisible
iff $M$ is $\mathcal{U}$-divisible.
\end{ex}

\medskip

According to the definition above, the perpendicular category $\mathcal{X}_\mathcal{U}$ of a class of bound modules $\mathcal{U}$ consists of the $\mathcal{U}$-torsion-free and $\mathcal{U}$-divisible modules.
It can  also be interpreted as the category of modules over the universal localization of $R$ at $\mathcal{U}$, as shown by
Crawley-Boevey \cite[Property~2.5]{CrawleyBoeveytubes} in a
slightly less general situation.

\begin{prop}\label{prop:Ralphamodules}
Let $R$ be a ring. Let $\mathcal{U}$ be a class of bound right
$R$-modules. The following statements are equivalent for
$M\in\rMod R.$
\begin{enumerate}[(1)]
\item $M\in\rMod R_{\mathcal{U}}$. \item $1_M\otimes_R\alpha_U^*$
is invertible for all morphisms $\alpha_U\in{\mathcal{U}}$. \item
$\Tor_R^1(M,\Tr U)=M\otimes_R\Tr U=0$ for every right $R$-module
$U\in\mathcal{U}.$ \item $\Hom_R(U,M)=\Ext_R^1(U,M)=0$ for every
right $R$-module $U\in{\mathcal{U}}$.
\end{enumerate}
\end{prop}

\begin{pro}
For the implication $(1)\Leftrightarrow (2)$, see the proof of
\cite[Theorem~4.7]{Schofieldbook}.

\noindent $(2)\Leftrightarrow (3)\Leftrightarrow  (4):$ Take
$\alpha_U\in{\mathcal{U}}.$ Then from
\begin{eqnarray*}
& \exseqmap PQ{U}{\alpha_U} &\\
& \exseqmap {Q^*}{P^*}{\Tr U}{\alpha_U^*}, &
\end{eqnarray*}
we get the following commutative diagram of exact sequences
$$\xymatrix@C=0.3cm@R=0.4cm{ 0\ar[r] & \Tor_1^R(M,\Tr U)\ar[r] &
M\otimes_R Q^*\ar[r]^{1_M\otimes_R \alpha_U^*}\ar[d] & M\otimes_R
P^*\ar[r]\ar[d] &
M\otimes_R\Tr U\ar[r] & 0 \\
0\ar[r] & \Hom_R(U,M)\ar[r] & \Hom_R(Q,M)\ar[r] &
\Hom_R(P,M)\ar[r] & \Ext_R^1(U,M)\ar[r] & 0}$$ where the vertical
arrows are isomorphisms. Hence $1_M\otimes_R \alpha_U^*$ is an
isomorphism if and only if $\Tor_R^1(M,\Tr U)=M\otimes_R\Tr U=0$
if and only if $\Hom_R(U,M)=\Ext_R^1(U,M)=0$.
\end{pro}

\begin{prop}[{cf.
\cite[Remark~5.8]{Angeleriherberatrlifaj2}}] Let
$\mathcal{U}$ be a class of bound right $R$-modules, and
$\Tr\mathcal{U}=\{ \Tr U \mid U\in\mathcal{U}\}$. Then the
class $\mathcal{U}^\perp$  of  $\mathcal{U}$-divisible
modules
 is a tilting class.
 Moreover,  the class of
$\Tr\mathcal{U}$-torsion-free
modules
 is a
cotilting class of  left $R$-modules. More precisely, it is the cotilting
class of cofinite type that corresponds to $\mathcal{U}^\perp$ under the bijective correspondence
 from  \cite[Theorem~2.2]{Angeleriherberatrlifaj1}.
\end{prop}
\begin{pro}
For the first statement,  see for example
\cite[Corollary~5.1.16]{GobelTrlifaj}. For the second
statement, recall that the correspondence in
\cite[Theorem~2.2]{Angeleriherberatrlifaj1} sends the
tilting class $\mathcal{U}^\perp$ to the cotilting class
$\mathcal{U}^\intercal =\{_RX\,\mid\, \Tor_1^R(U,X)=0 \
\text{for all}\, U\in\mathcal{U}\}.$

If $U\in\mathcal{U},$ and $\exseqmap PQU\alpha$ is a projective
presentation of $U$ with $P$ and $Q$ finitely presented, we obtain
the exact sequences $\exseqmap{Q^*}{P^*}{\Tr U}{}$,  and
$$0\longrightarrow {\Tor_1^R(U,X)}\longrightarrow{P\otimes X}\longrightarrow{Q\otimes X}\longrightarrow{U\otimes X}\longrightarrow 0$$
$$0\to {\Hom_R(\Tr U,X)}\to{\Hom_R(P^*,X)}\to{\Hom_R(Q^*,X)}\to{\Ext_R^1(\Tr U,X)}\to 0.$$
Since $P\otimes_RX\cong\Hom_R(P^*,X)$ and $Q\otimes_RX\cong
\Hom_R(Q^*,X)$ are naturally isomorphic, we get
$\Tor_1^R(U,X)\cong \Hom_R(\Tr U,X).$ Therefore
$X\in{}\mathcal{U}^\intercal$ if and only if $\Hom_R(\Tr U,X)=0$
for all $U\in\mathcal{U},$ that is, $X$ is
$\Tr\mathcal{U}$-torsion-free.
\end{pro}

\medskip

When $\lambda:R\to R_{\mathcal{U}}$ is injective and
pd$R_{\mathcal{U}}\le 1$, we have a tilting module
$R_{\mathcal{U}}\oplus R_{\mathcal{U}}/R$ by Theorem
\ref{theo:tiltingepimorphism}. In general, however, its tilting
class $\Gen R_{\mathcal{U}}$ does not coincide with the tilting
class ${\mathcal{U}}^\perp$, as we will see in Example
\ref{ex:tame}.
 The next result describes the case when $\Gen R_{\mathcal{U}}={\mathcal{U}}^\perp$. We first need some preliminaries.

 \begin{des}

An ascending chain $(N_\nu|\nu<\kappa)$ of submodules of a right
$R$-module $N$ indexed by a cardinal $\kappa$ is called
\emph{continuous} if
$N_\nu=\operatornamewithlimits{\cup}\limits_{\beta<\nu}N_\beta$
for all limit ordinals $\nu<\kappa.$ The continuous  chain is
called a \emph{filtration} of $N$ if $N_0=0$ and
$N=\operatornamewithlimits{\cup}\limits_{\nu<\kappa} N_\nu.$

Given a class $\mathcal{U}$ of right $R$-modules, we say that a
right $R$-module $N$ is \emph{$\mathcal{U}$-filtered} if it admits
a filtration $(N_\nu|\nu<\kappa)$ such that $N_{\nu+1}/N_\nu$ is
isomorphic to some module in $\mathcal{U}$ for every $\nu<\kappa.$
\end{des}

The following result is well known.

\begin{lem}[{\cite[Lemma~3.1.2]{GobelTrlifaj}}]\label{lem:ekloflemma}
Let $M$ be a right $R$-module, and let $\mathcal{U}$ be a class of right $R$-modules
such that $M\in\mathcal{U}^\perp$.
If $N$ is a {$\mathcal{U}$-filtered} right $R$-module, then $M\in N^\perp$.
\end{lem}

\begin{theo}\label{theo:envelope}
Let $R$ be a ring. Let $\mathcal{U}$ be a class of   bound right
$R$-modules. Let further $\lambda:R\to S$ be an injective ring
epimorphim with $\Tor_1^R(S,S)=0$ and $\pd S_R\le 1$. The
following statements are equivalent.
\begin{enumerate}[(1)]
\item $\Gen S_R=\mathcal{U}^\perp$.
\item The   map $\lambda\colon R\rightarrow S$ is a $\mathcal{U}^\perp$-(pre)envelope.
\item $S_R\in \mathcal{U}^\perp$, and every (pure-injective) module $M\in\mathcal{U}^\perp$ belongs to
$(S/R)^\perp$.
\end{enumerate}
In particular, conditions $(1)-(3)$ hold true if $S_R\in \mathcal{U}^\perp$ and $S/R$ is a direct limit of $\mathcal{U}$-filtered right $R$-modules.
\end{theo}
\begin{pro}
We already know by
Theorem~\ref{theo:tiltingepimorphism} that
$T=S\oplus S/R$ is a tilting right $R$-module with
$\Gen T=\Gen S_R=(S/R)^\perp$.

\noindent $(1)\Rightarrow (2):$ If $M\in {\mathcal{U}}^\perp=\Gen S_R$, then
$\Ext_R^1(S/R,M)=0$, and therefore  $\Hom_R(\lambda, M)$ is surjective. So, $\lambda: R\to S$ is a ${\mathcal
U}^\perp$-preenvelope. Suppose now that
$g\in\End_R(S)$ satisfies
$\lambda=g\lambda.$ Since $\rMod S$ is a
full subcategory of $\rMod R$,
$g\in \End_{S}(S).$ Now,
since $g(1)=1,$ we get that $g$ is the identity and therefore an
isomorphism. So $\lambda$ is even a ${\mathcal
U}^\perp$-envelope.

\noindent $(2)\Rightarrow (1):$ By the definition of a
preenvelope, we have that $S_R$ belongs to $\mathcal{U}^\perp.$
Since $\mathcal{U}^\perp$ is a torsion class, it follows $\Gen
S_R\subseteq \mathcal{U}^\perp.$ For the reverse inclusion, note that $S_R$ is a generator of
$\mathcal{U}^\perp$ by \cite[Lemma~1.1]{ATT}.

\noindent $(1)\Rightarrow(3)$ follows from $(S/R)^\perp=\Gen S_R$.

\noindent $(3)\Rightarrow (1):$ Since $\mathcal{U}^\perp$ is a
torsion class, we deduce as above $\Gen S_R\subseteq
\mathcal{U}^\perp.$ To prove equality, note   that both classes
are tilting classes. By \cite[Theorem~1.6]{BazzoniHerbera}  it
follows that they are both \emph{definable classes}, that is, they
are closed under direct products, direct limits, and pure
submodules. As noted in Section~2 of \cite{BazzoniHerbera}, a
combination of Ziegler's result \cite[Theorem~6.9]{Ziegler} and
Keisler-Shelah Theorem (cf. \cite{Keisler} and \cite{Shelah})
implies that $\Gen S_R$ and  $\mathcal{U}^\perp$ coincide if and
only if they contain the same pure-injective right $R$-modules.
The latter holds true by (3).

\smallskip

We now prove the last statement. Suppose that $S_R\in \mathcal{U}^\perp$ and $S/R=\lim\limits_{\longrightarrow} N_i$ where all
$N_i$ are $\mathcal{U}$-filtered right $R$-modules. By condition (3)  it is enough
to show that every pure-injective module $M\in\mathcal{U}^\perp$
belongs to $(S/R)^\perp.$
Now, for such module  $M$ we have
$$\Ext_R^1(S/R,M)=\Ext_R^1(\lim\limits_{\longrightarrow} N_i,M)\cong
\lim\limits_{\longleftarrow}\Ext_R^1(N_i,M).$$ Since all
$N_i$ are $\mathcal{U}$-filtered,
$\Ext_R^1(N_i,M)=0$ by Lemma~\ref{lem:ekloflemma}, thus
$\Ext_R^1(S/R,M)=0.$
 \end{pro}

\begin{coro}\label{theo:tiltingwithclass}
Let $R$ be a ring. Let $\mathcal{U}$ be a class of   bound right
$R$-modules. Suppose that $R$ embeds in $R_{\mathcal{U}},$ and
$\pd (R_{\mathcal{U}})_R\leq 1$. Assume further that
$R_{\mathcal{U}}/R$ is a direct limit of $\mathcal{U}$-filtered right $R$-modules. Then
$T_{\mathcal{U}}=R_{\mathcal{U}}\oplus R_{\mathcal{U}}/R$ is a
tilting right $R$-module with $\Gen T_{\mathcal{U}}=\Gen
(R_{\mathcal{U}})_R=\mathcal{U}^\perp.$
\end{coro}
\begin{pro}
Notice that
$R_{\mathcal{U}}\in\mathcal{U}^\perp$ because $R_{\mathcal{U}}$ is
a right $R_{\mathcal{U}}$-module, see
Proposition~\ref{prop:Ralphamodules}. So, the statement follows immediately from Theorem \ref{theo:envelope} and Remark \ref{rems:universalisepi}(3).
\end{pro}

We now give some applications of the last results.
We start with an extension of
\cite[Proposition~6.4]{Angeleriherberatrlifaj2}.

\begin{coro}\label{coro:Oretiltingclass}
Let $R$ be a ring. Let $\mathfrak{U}$ be a left Ore set of
non-zero-divisors of $R.$ Then $\pd(\mathfrak{U}^{-1}R_R)\leq 1$
if and only if $\Gen(\mathfrak{U}^{-1}R_R)$ coincides with the
class of $\mathfrak{U}$-divisible right $R$-modules. In this case
$T_{\mathfrak{U}}=\mathfrak{U}^{-1}R\oplus \mathfrak{U}^{-1}R/R$
is a tilting right $R$-module whose tilting class coincides with
the class of $\mathfrak{U}$-divisible right $R$-modules.
\end{coro}

\begin{pro}
Suppose $\pd(\mathfrak{U}^{-1}R_R)\leq 1.$ Since $\mathfrak{U}$
consists of non-zero-divisors, $R$ embeds in $\mathfrak{U^{-1}}R.$
Setting $\mathcal{U}=\{R/uR\mid u\in\mathfrak{U}\},$ we know by
Example~\ref{ex:oreisuniversal} that $\mathcal{U}$ is a set of bound
right
 $R$-modules and $R_\mathcal{U}\cong\mathfrak{U}^{-1}R.$ On
the other hand, given $u,v\in\mathfrak{U},$ there exist
$z\in\mathfrak{U},\ w\in R,$ such that $wu=zv.$ Then
$$u^{-1}R+v^{-1}R\subseteq (zv)^{-1}R=(wu)^{-1}R.$$ Hence every
finitely generated right submodule of $\mathfrak{U}^{-1}R$ is
contained in $u^{-1}R$ for some $u\in\mathfrak{U}.$ Therefore
$\mathfrak{U}^{-1}R/R=\operatornamewithlimits{\lim}\limits_{\substack{\longrightarrow\\
u\in\mathfrak{U}}} u^{-1}R/R.$ Moreover, notice that for every
$u\in\mathfrak{U},$ $u^{-1}R/R\cong R/uR,$ thus,
$\mathfrak{U}^{-1}R/R$ is a direct limit of the
$\mathcal{U}$-filtered modules $u^{-1}R/R.$ Then, applying
Corollary~\ref{theo:tiltingwithclass}, we obtain that
$T_{\mathfrak{U}}=\mathfrak{U}^{-1}R\oplus \mathfrak{U}^{-1}R/R$
is a tilting right $R$-module and  $\Gen (T_{\mathfrak{U}})_R=\Gen
(\mathfrak{U}^{-1}R)_R=\mathcal{U}^\perp.$

The proof of the other implication is given in the last paragraph
of the proof of  \cite[Proposition~6.4]{Angeleriherberatrlifaj2}.
\end{pro}

\begin{rem}
If $\mathfrak{U}$ is a twosided Ore set of
non-zero-divisors, then $\pd({}_R\mathfrak{U}^{-1}R)\leq 1$ if and
only if $\Gen({}_R\mathfrak{U}^{-1}R)$ coincides with the class of
$\mathfrak{U}$-divisible left $R$-modules. In fact, in this
case $\mathfrak{U}^{-1}R=R\mathfrak{U}^{-1}=R_\mathcal{U},$ and we can apply
 the left version of Corollary~\ref{theo:tiltingwithclass} on
${}_RR\mathfrak{U}^{-1}.$

However, if  $\mathfrak{U}$ is just a left Ore set of
non-zero-divisors of $R,$
and
$\pd({}_R\mathfrak{U}^{-1}R)\leq 1,$ then
$T_{\mathfrak{U}}=\mathfrak{U}^{-1}R\oplus \mathfrak{U}^{-1}R/R$
is a  tilting left $R$-module by Theorem~\ref{theo:tiltingepimorphism}, but we cannot compute
$T_{\mathfrak{U}}^\perp$ as we don't know whether
$\mathfrak{U}^{-1}R/R$ can be written as a direct limit of $\{R/Ru\mid
u\in\mathfrak{U}\}$-filtered left $R$-modules.

Stronger results will be obtained in
Theorem~\ref{theo:tiltingprimenoetherianrings} under the assumption that $R$ is a hereditary noetherian prime ring.
\end{rem}

\smallskip

\begin{coro}[{\cite[Remark~6.3.17]{GobelTrlifaj}}] \label{coro:alltiltingvaluation}
Let $R$ be a commutative valuation domain with field of fractions
$Q.$ Suppose that $\pd (R_\mathfrak{p})_R\leq 1$  for each prime
ideal $\mathfrak{p}$ of $R$ (equivalently, suppose that
$R_\mathfrak{p}$ is countably generated as an $R$-module for every
prime ideal $\mathfrak{p}$ of $R$).   Then the set
$\mathbb{T}=\{T_\mathfrak{p}=R_\mathfrak{p}\oplus
R_\mathfrak{p}/R\mid \mathfrak{p} \in \Spec(R) \}$ is a
representative set up to equivalence of the class of all tilting
$R$-modules.
\end{coro}

\begin{pro}
For each  prime ideal $\mathfrak{p}$ of $R$, let
$\mathfrak{U}_\mathfrak{p}=R\setminus\mathfrak{p}.$ By
Corollary~\ref{coro:Oretiltingclass} we know that
$T_\mathfrak{p}=R_\mathfrak{p}\oplus R_\mathfrak{p}/R$ is a
tilting $R$-module and $T_\mathfrak{p}^\perp$ equals  the class of
$\mathfrak{U}_\mathfrak{p}$-divisible $R$-modules.

It is known that the set of  Fuchs tilting modules
$\{\delta_{\mathfrak{U}_\mathfrak{p}}\mid \mathfrak{p}\in\textrm{
Spec($R$)}\}$ is a representative set up to equivalence of the
class of all tilting $R$-modules, and
$\delta_{\mathfrak{U}_\mathfrak{p}}^\perp$ is the class of
$\mathfrak{U}_\mathfrak{p}$-divisible $R$-modules
\cite[Theorem~6.2.21]{GobelTrlifaj}.

The assumption that $\pd (R_\mathfrak{p})_R\leq 1$  for each prime
ideal $\mathfrak{p}$ of $R$ is satisfied if and only if
$R_\mathfrak{p}$ is countably generated as an $R$-module for every
prime ideal $\mathfrak{p}$ of $R.$ In  fact, if $R$ is a
commutative local ring and $\mathfrak{U}$ is a multiplicative
subset of non-zero-divisors, then $\pd R\mathfrak{U}^{-1}_R\leq 1$
if and only if $R\mathfrak{U}^{-1}_R$ is countably generated
$R$-module \cite[Page~531]{Angeleriherberatrlifaj2}. In the
particular  case when $R$ is a valuation domain see
\cite[Theorem~IV.3.1]{FuchsSalcebook}
\end{pro}

\smallskip

\begin{ex}\label{ex:freegroupalgebra}
(see also \cite[Example~0.2]{NeemanNoncommutativeII})
Let $X$ be a nonempty set. Let $G$ be the free group on $X.$ Let
$k$ be a field. Consider the free algebra $R=\freealgebra kX$ and
the free group algebra $kG$ with the natural embedding
$\freealgebra kX\hookrightarrow kG$ which sends  $x\mapsto x$ for
every $x\in X.$ Let $\mathcal{X}=\{\freealgebra kX/x\freealgebra
kX \mid x\in X\}.$ Then $T_X= kG\oplus kG/\freealgebra kX$ is a
tilting right  $R$-module with $T_X^\perp=\mathcal{X}^\perp$.

In fact, if $\Sigma=\{\alpha_x\mid x\in X\}$ where $\alpha_x\colon
\freealgebra kX\rightarrow \freealgebra kX$ is defined by
$p\mapsto xp$, then $kG$ can be
regarded as the universal localization of $R$ at $\mathcal{X}.$
Since $\freealgebra kX$ is hereditary, $\pd(kG)\leq 1,$ so $T_X$
is a tilting right $R$-module by
Remark~\ref{rems:universalisepi}(3).

 Note that
$\mathcal{X}$ is a set of bound right   $\freealgebra kX$-modules.
 Since $kG$ is a right $kG$-module,
$kG\in\mathcal{X}^\perp$ by Proposition~\ref{prop:Ralphamodules}.
We now verify condition (2) in Theorem~\ref{theo:envelope}. Let
$M\in\mathcal{X}^\perp.$ We have to show that for every
$\freealgebra kX\stackrel{\widetilde{f}}{\rightarrow}M,$
 there exists $f\colon kG\rightarrow M$ extending
$\widetilde{f}.$ We will define $f$ on the elements of $G$ and
then extend it by linearity.

Every element $g\in G$ can be uniquely expressed as a word of the
form
\begin{equation}\label{eq:normalformfreegroup}
g=x_1^{e_1}\dotsb x_r^{e_r} \textrm{ where } x_i\in X,\ e_i=\pm 1
\textrm{ and } x_i\neq x_{i+1} \textrm{ if } e_i=-e_{i+1}.
\end{equation}
We proceed by induction on the length of $g.$ If $r=0,$ that is,
$g=1,$ then we define $f(1)$ as $\widetilde{f}(1).$ Let $r+1>0$
and suppose we have defined $f(g)$ for all $g\in G$ of length
$\leq r.$ Let $g$ be a word of length $r+1,$ suppose
$g=x^{e_1}\dotsb x_{r+1}^{e_{r+1}}=hx_{r+1}^{e_{r+1}}$ as in
\eqref{eq:normalformfreegroup}. Notice that for each $m\in M$ and
$x\in X$ there exists $n\in M$ such that $m=nx.$  This can be seen
applying $\Hom_R({}_-,M)$ to the short exact sequence determined
by $\alpha_x$ and the fact that $M\in\mathcal{X}^\perp.$ So fix
$n\in M$ such that $nx_{r+1}=f(h).$ Then
$$f(g)=\left\{\begin{array}{ll}
f(h)x^{e_{r+1}}_{r+1} & \textrm{ if } e_{r+1}=1 \\
n  & \textrm{ if } e_{r+1}=-1.
\end{array}\right.$$
Hence $\freealgebra kX\hookrightarrow kG$ is an  $\mathcal{X}^\perp$-preenvelope, and $T_X^\perp=\Gen kG=\mathcal{X}^\perp$.

Finally, observe that $kG$ is not a flat right (left)
$\freealgebra kX$-module if $|X|\geq 2$ . Indeed, let $x\neq y\in
X.$ Consider the unique embedding of left (right) $\freealgebra
kX$-modules such that
\begin{equation*}
\begin{array}{ccc} \freealgebra kX \oplus \freealgebra kX & \stackrel{\alpha}{\longrightarrow} & \freealgebra kX \\
(1,0) & \longmapsto & x\\
(0,1) & \longmapsto & y
\end{array}
\end{equation*}

Consider $1_{kG}\otimes \alpha\colon kG\oplus kG\longrightarrow
kG.$ Then $(x^{-1},0)$  and $(0,y^{-1})$ have the same image $1.$
Thus $1_{kG}\otimes \alpha$ ($\alpha\otimes 1_{kG}$) is not
injective.

It can be seen that $T_X$ is also a tilting left $R$-module with
 ${}_RT_X^\perp=\{\freealgebra
kX/\freealgebra kX x\mid x\in X\}^\perp.$
\end{ex}

\section{Projective rank functions}

We recall some notions and results from \cite{Schofieldbook} and \cite{CrawleyBoeveytubes}. For details, we refer to
\cite[Theorems~1.11,~1.16,~1.18,~1.22,~5.1,~5.5]{Schofieldbook}  and
\cite[Theorem~1.4]{CrawleyBoeveytubes}.

\begin{des}
 Let $R$ be a ring. We denote by
$K_0(R)$ the \emph{Grothendieck group} of finitely generated
projective right $R$-modules modulo direct sums, that is, the
abelian group generated by the isomorphism classes $[P]$ of
$P\in\mathcal{P}_R$ modulo the relations $[P]+[Q]-[P\oplus Q]$
for all $P,Q\in \mathcal{P}_R.$

\smallskip

(a) A \emph{(projective) rank function} on a ring $R$ is a morphism of
groups $\rho\colon K_0(R)\rightarrow \mathbb{R}$ such
that\begin{enumerate}[(i)] \item $\rho([P])\geq 0$ for all
$P\in\mathcal{P}_R$, \item $\rho([R])=1.$
\end{enumerate}
If $\rho([P])> 0$ for every nonzero $P\in\mathcal{P}_R,$ we say that
$\rho$ is a \emph{faithful rank function}. For  sake of
simplicity we will write $\rho(P)$ instead of $\rho([P])$.

\smallskip

(b) Let $\alpha\colon P\rightarrow Q$ be a morphism between finitely
generated projective right $R$-modules. Consider the finitely
generated projective right $R$-modules $P'$ such that there exist
morphisms $\beta,\gamma$ making the following diagram commutative
\begin{eqnarray}\label{eq:rankofmaps}
\xymatrix@C=0.3cm@R=0.3cm{ P\ar[rr]^\alpha\ar[dr]_\beta & & Q\\ &
P'\ar[ur]_\gamma }
\end{eqnarray}We define the \emph{inner rank} of $\alpha$ as
$\rho(\alpha)=\inf\{\rho(P')\mid P' \textrm{ satisfies
\eqref{eq:rankofmaps}}\}.$

\smallskip

(c) We say that a morphism between finitely generated projective
right $R$-modules $\alpha\colon P\rightarrow Q$ is \emph{full} in
case $\rho(\alpha)=\rho(P)=\rho(Q).$ We denote the localization of
$R$ at the set of all full morphisms by $R_\rho,$ and we  call it
the \emph{universal localization of $R$ at $\rho$}. If $\rho$ is
faithful and $\alpha$ is full, we define $\alpha$ to be an
\emph{atomic full morphism} if, in any nontrivial factorization as
in \eqref{eq:rankofmaps}, we have
$\rho(P')>\rho(P)=\rho(Q)=\rho(\alpha).$
\end{des}

\begin{theo}\label{theo:extensionandkernel}
Let $R$ be a hereditary ring with a faithful rank function $\rho.$ Let
$\Sigma$ be a collection of full maps. Then  $R$ embeds in $R_\Sigma$.
\end{theo}

\begin{de}\label{de:exact}
Suppose $R$ is a semihereditary ring with a faithful rank function
$\rho.$ Notice  that every full map is injective under these
assumptions since $P\cong\ker\gamma\oplus\im\gamma$ for every
morphism $\gamma\colon P\rightarrow Q$ between finitely generated
projective right $R$-modules. Let $M$ be a right $R$-module. We
say that $M$ is \emph{$\rho$-torsion} if $M$ is the cokernel of a
full morphism. We say that $M$ is \emph{$\rho$-simple} if $M$ is
the cokernel of an atomic full morphism.

Of course, the
$\rho$-torsion (and the $\rho$-simple) modules are bound right
$R$-modules. Moreover, if $R$ is hereditary, then
the $\rho$-torsion modules form an exact abelian length category (that is, every object has finite length) whose simple objects are the $\rho$-simples.

For a characterization of $\rho$-torsion, $\rho$-simple modules in
the hereditary case see \cite[Definition~1.3]{CrawleyBoeveytubes}.
\end{de}

\begin{theo}[{\cite{Schofieldquivers}, \cite[Theorem~12.6]{Schofieldbook}}]
\label{theo:universallocalizationwithfiltration}
 Let $R$ be a hereditary ring with a faithful rank function
$\rho$ such that $R_\rho$ is a simple artinian ring. Let
$\mathcal{U}$ be a class of $\rho$-simple modules. The following
statements hold true.
\begin{enumerate}
\item As a right $R$-module,
$R_\mathcal{U}/R$ is a directed union of finitely presented
modules $N_i$
 such that each $N_i$ is a finite extension of modules from
$\mathcal{U}.$ \item As a left $R$-module, $R_\mathcal{U}/R$ is a
directed union of finitely presented modules $M_j$ such that each
$M_j$ is a finite extension of modules of the form $\Tr U$ with
$U\in \mathcal{U}.$
\end{enumerate}
\end{theo}

In particular, the Theorem above applies in the following situation.

\begin{theo}\label{theo:universallocalizationsimpleartinian}
Let $R$ be a hereditary ring. Suppose $R$ has a unique rank
function $\rho,$ and suppose that $\rho$ takes values in
$\frac{1}{n}\mathbb{Z}$ for some positive integer $n.$ Then the
universal localization $R_\rho$ of $R$ at $\rho$ is a simple
artinian ring.
\end{theo}

\bigskip

Let us now apply the results above.
Combining  Theorems~\ref{theo:extensionandkernel}
and \ref{theo:universallocalizationwithfiltration} with   Remark \ref{rems:universalisepi}(3) and Corollary~\ref{theo:tiltingwithclass},
we immediately obtain the following.

\begin{coro}\label{coro:rhotorsiongivetilting} Let $R$ be a hereditary ring
with a faithful rank function $\rho.$ The following statements hold true.
\begin{enumerate}[(1)]
\item If $\mathcal{V}$ is a class of $\rho$-torsion right
$R$-modules, then $T_\mathcal{V}=R_\mathcal{V}\oplus
R_\mathcal{V}/R$ is a tilting right $R$-module. \item Suppose
$R_\rho$ is simple artinian.  If $\mathcal{U}$ consists of
$\rho$-simple modules, then $T_\mathcal{U}=R_\mathcal{U}\oplus
R_\mathcal{U}/R$ is a tilting right  $R$-module with tilting class
$T_\mathcal{U}^\perp=\mathcal{U}^\perp$, and it is also  a tilting
left  $R$-module with tilting class
$_R\,T_\mathcal{U}\,^\perp=(\Tr\mathcal{U}) ^\perp$.
\end{enumerate}
\end{coro}

In general, if $\mathcal{V}$ is just a class of $\rho$-torsion right $R$-modules,
the tilting class   $\Gen T_\mathcal{V}$  differs from $\mathcal{V}^\perp$, as we are going to see next.

\begin{ex}\label{ex:tame}
Let  $R$ be an (indecomposable) tame hereditary algebra.
Then $K_0(R)$ is the free
abelian group with basis the (finite number of) isomorphism
classes of simple right $R$-modules $\{S_1,\dotsc,S_n\},$ and there is a bilinear form $B_R\colon K_0(R)\times
K_0(R)\rightarrow \mathbb{Z}$ given by $$B_R([M],[N])=\dim_k
\Hom_R(M,N)-\dim_k \Ext_R^1(M,N)$$ with  corresponding
quadratic form $\chi_{{}_R}\colon
\mathbb{Q}\otimes_\mathbb{Z}K_0(R)\rightarrow \mathbb{Q}$.
Since $R$ is of tame representation type, $\chi_{{}_R}$ is positive semidefinite
but not positive definite.

Moreover, the $\mathbb{Q}$-subspace
$N\leq \mathbb{Q}\otimes_{\mathbb{Z}} K_0(R)$ formed by the
radical vectors of $B_R$ is one-dimensional and can be generated
by a vector $v$ with coordinates $(v_1,\dotsc,v_n)\in\mathbb{N}^n$
in the basis $\{[S_1],\dotsc,[S_n]\}$ with at least one component
$v_i=1,$ see \cite{RingelTamealgebras}.
Hence $\chi_{{}_R}(v)=0$ and any other $w$ such that
$\chi_{{}_R}(w)=0$ is a $\mathbb{Q}$-multiple of $v.$

Following \cite[Section~4]{CrawleyBoeveytubes}, we define a
faithful rank function $\partial_R\colon K_0(R)\rightarrow
\mathbb{Q}$  by
$$\partial_R([M])=\frac{B_R([M],v)}{B_R([R],v)}.$$ It is called
the \emph{normalized defect for $R$}.

The indecomposable $\partial_R$-torsion modules coincide with the
regular modules and the $\partial_R$-simple modules with the
{simple regular modules}. Moreover, it was  shown in
\cite[Lemma~4.4]{CrawleyBoeveytubes} that $R_{\partial_R}$ is a
simple artinian ring. From Corollary
\ref{coro:rhotorsiongivetilting}(2) we  infer that
 for every set
$\mathcal{U}$ of  simple regular right $R$-modules, the right
$R$-module $T_\mathcal{U}=R_\mathcal{U}\oplus R_\mathcal{U}/R$ is
a tilting module with tilting class
$T_\mathcal{U}^\perp=\mathcal{U}^\perp.$

Let us now assume  that there is a stable tube $\mathbf{t}_{\nu}$ of width at
least $3.$ Let $S$ be a simple regular module such that $[S]\in
\mathbf{t}_{\nu}.$ Consider the modules $\tau S,\ \tau^{-}S,\ S[2]$ and
$\tau S[2].$

The $\partial_R$-torsion module $S[2]$ is an
extension of the $\partial_R$-simple modules $S$ and $\tau^-S$.
Similarly, $\tau S[2]$ is an
extension of the $\partial_R$-simple modules $\tau S$ and $S$.
Hence we can suppose that $S[2],\ \tau S[2]$
have finite projective presentations \begin{eqnarray*} \exseqmap
{P_1} {Q_1} {\tau S[2]}{\alpha}\\ \exseqmap {P_2} {Q_2}
{S[2]}{\beta}
\end{eqnarray*}
where $\alpha$ and $\beta$ are full morphisms, and that
$\alpha=\delta\gamma,$  $\beta=\varepsilon\delta$ where
$\gamma,\delta,\varepsilon$ are full atomic morphisms with
$\coker\gamma= \tau S,\ \coker\delta=S$ and
$\coker\varepsilon=\tau^-S.$

Let $\mathcal{V}=\{S[2],\ \tau S[2]\}$ and $\mathcal{U}=\{S,\ \tau
S,\ \tau^- S\}.$ Then $R_\mathcal{V}=R_{\{\alpha,\beta\}}\cong
R_{\{\gamma,\delta,\varepsilon\}}=R_\mathcal{U}.$ Therefore
$T_\mathcal{V}= R_\mathcal{V}\oplus R_\mathcal{V}/R$ is a tilting
right $R$-module with $T_\mathcal{V}^\perp=\mathcal{U}^\perp$. But
$\mathcal{U}^\perp$ is different from $\mathcal{V}^\perp.$ In fact, using
the AR-formula we have
\begin{eqnarray*} \Ext_R^1(S[2],S) & \cong & D\Hom_R(S,\tau S[2])=0 \\
\Ext_R^1(\tau S[2],S) & \cong & D\Hom_R(S,\tau^2S[2])=0.
\end{eqnarray*}
Hence $S\in\mathcal{V}^\perp.$ But \begin{eqnarray*}
\Ext_R^1(\tau^-S,S)\cong D\Hom_R(S,\tau\tau^-S)\cong
D\Hom_R(S,S)\neq 0,
\end{eqnarray*}
that is, $S\notin \mathcal{U}^\perp.$
\end{ex}

\begin{rem}
Many rings  satisfy the conditions in
Corollary~\ref{coro:rhotorsiongivetilting}(2), for example
tame hereditary algebras (see above),  Dedekind prime
rings~\cite[Proposition~3.1~(2)]{CrawleyBoeveytubes} and
firs. Hereditary local rings, free algebras and free group
algebras are examples of firs.
\end{rem}


\section{Noetherian prime rings}

\begin{nota}
Let $R$ be a right order in a semisimple  ring $A,$ i.e.
$A$ is the right Ore localization of $R$ at the set
$\mathcal{C}_R$ of regular elements of $R.$ Let $n$ be the
length of $A$ as a right $A$-module. We define the rank
function $u\colon K_0(R)\rightarrow \frac{1}{n}\mathbb{Z}$
given by
$$u(P)=\frac{\leng(P\otimes_R A_A)}{n},$$ called
the \emph{normalized uniform dimension}, see \cite{CrawleyBoeveytubes}.

Let $\mathcal{U}_r$ be a set of representatives of all isomorphism
classes of finitely presented simple right $R$-modules. Let
$\mathcal{V}_r$ be a set of representatives of all isomorphism
classes of finitely presented torsion right $R$-modules. Let
finally $\mathcal{D}_r=\{R/sR\mid s\in\mathcal{C}_R\}.$ In the
same way we define $\mathcal{U}_l,\ \mathcal{V}_l,\
\mathcal{D}_l.$
\end{nota}

\begin{rems}
\begin{enumerate}[(1)]
\item Since $A$ is  the right Ore localization of $R$ at
$\mathcal{C}_R,$ the left $R$-module ${}_RA$ is  flat. Moreover a
right $R$-module $V_R$ is torsion if and only if $V\otimes_R A=0.$

\item Since every projective right $R$-module is torsionfree, $u$
is a faithful rank function.
\end{enumerate}
\end{rems}

Recall the following result which can be found, for example, in
\cite[Corollary~2.2.12]{Jategaonkarbook}.

\begin{lem}\label{lem:simplearetorsion}
\label{lem:rightordersandsimplemodules} Let $R$ be a right order
in a simple artinian ring. If there exists a simple torsionfree
right $R$-module, then $R$ itself is a simple artinian ring.
\end{lem}

The following result is a generalization of the foregoing Lemma to
the semisimple situation.

\begin{lem}\label{lem:semisimplerightordersandsimplemodules}
Let $R$ be a right order in a semisimple ring $A.$ If there exists
a simple torsionfree right $R$-module, then there exists a
primitive central idempotent $e$ of $A$ such that $eRe$ is a
simple artinian ring.
\end{lem}
\begin{pro}
Suppose $M$ is a simple torsionfree right $R$-module. It is known
that there exists a right ideal $I$ of $R$ such that $M$ and $I$
have isomorphic essential submodules, see for example
\cite[Proposition~2.2.11]{Jategaonkarbook}. Hence, since $M$ is
simple, $M$ embeds in $R.$ So we can suppose $M$ is a right ideal
of $R.$ There exists a primitive central idempotent $e$ of $A$
such that $Me\neq 0.$ Then $Me\cong M$ as right $R$-modules, and
$Me$ is a simple torsion free right $eRe$-module. Notice that $eRe$ is
a right order in the simple artinian ring $eAe.$ Thus  $eRe$ is
simple artinian by Lemma~\ref{lem:rightordersandsimplemodules}.
\end{pro}

\begin{prop}[{cf. \cite[Section~3]{CrawleyBoeveytubes}}] \label{theo:semihereditarytilting}
Let $R$ be a semihereditary right order in a semisimple ring $A$.
Suppose there is no primitive central idempotent $e$ of $A$ such
that $eRe$ is simple artinian.   Then
\begin{enumerate}[(1)] \item The class of
finitely presented torsion right $R$-modules coincides with
the class of $u$-torsion modules. \item The class of
finitely presented simple right $R$-modules coincides with
the class of $u$-simple modules.
\item $A$ equals $R_u$, the universal localization of $R$ at $u$.
\end{enumerate}
\end{prop}

\begin{pro}
(1) Given a finitely presented torsion right $R$-module $V_R$
(hence $\pd V_R=1$) with finite projective presentation $\exseqmap
PQV\alpha,$ applying ${}_-\otimes_R A,$ we get $\lexseqmap
{P\otimes_RA}{Q\otimes_R A}{V\otimes_RA=0}{\alpha\otimes 1_A}.$
Hence $u(P)=u(Q)=u(\alpha),$ and $V$ is a $u$-torsion module.
Conversely, if $V$ is a $u$-torsion module, then $V$ is a finitely
presented module with $\pd V_R=1.$ Let $\exseqmap PQV\alpha$ be a
 projective presentation of $V$ with $P$ and $Q$ finitely generated.
   Notice that
$\length(P\otimes_R A)=\length(Q\otimes_R A).$  Hence
$\alpha\otimes 1_A$ is an isomorphism and $V\otimes_R A=0.$ Thus
$V$ is a torsion right $R$-module.

(2) Let $U$ be a finitely presented simple right $R$-module  with
finite projective presentation $\exseqmap PQU\alpha.$ Since $U$ is
simple and there is no primitive central idempotent $e$ of $A$
such that $eRe$ is simple artinian,
Lemma~\ref{lem:semisimplerightordersandsimplemodules} implies $U$
is a torsion right $R$-module, and therefore $u$-torsion  with
$\pd U_R=1$ by (1). Now suppose
$$\xymatrix@C=0.3cm@R=0.3cm{ P\ar[rr]^\alpha\ar[dr]_\beta & & Q\\
& P'\ar[ur]_\gamma }$$ with $u(P')=u(P)=u(Q)=u(\alpha).$  Hence
$\length(Q\otimes_R A)=\length(P'\otimes_R A).$ Since
$\alpha\otimes 1_A$ is surjective, we get $\gamma\otimes 1_A$  is
an isomorphism. Hence we have the commutative diagram
$$\xymatrix@C=0.5cm@R=0.4cm{ P'\otimes_RA\ar[r]^\cong &  Q\otimes_RA \\
 P'\ar[r]_\gamma\ar[u] & Q\ar[u] }$$
where the vertical arrows are injective. Hence $\gamma$ is
injective. Clearly $\beta$ is injective. Now, since $U\cong Q/P$
is simple, we get that $\beta$ or $\gamma$ is an isomorphism. This shows that $U$ is $u$-simple.

On the other hand, if $U$ is a $u$-torsion module which is not a
simple module, then it contains a finitely generated submodule
$0\neq V\lvertneqq U.$ Suppose $\exseqmap PQU\alpha$ is a
projective presentation of $U$ with $P$ and $Q$ finitely
generated. Then there exists a finitely generated  submodule $0\ne
P'\lneqq Q$ such that $P'/P\cong V.$ Since $R$ is semihereditary,
$P'$ is a projective right $R$-module. Now  $\alpha$ factors
through $P'$ in the following
way $$\xymatrix@C=0.3cm@R=0.3cm{ P\ar[rr]^\alpha\ar@{^{(}->}[dr] & & Q\\
& P'\ar@{^{(}->}[ur] }$$   and $U$ cannot be a $u$-simple module
since $u(P)=u(Q)=u(P')$.

(3) If $\exseqmap PQW\alpha$ is an exact sequence with $W$
torsion and $P,Q$ finitely generated projective right
$R$-modules, then $\alpha\otimes 1_A\colon
P\otimes_RA\rightarrow Q\otimes_RA$ is an isomorphism.
Therefore condition (i) in
Theorem~\ref{def:universallocalization}  is satisfied.

Let $b\colon R\rightarrow B$ be a morphism of rings such that
$\alpha\otimes 1_B$ is invertible for every full morphism
$P\stackrel{\alpha}{\rightarrow} Q$. By (1), $\alpha$ is full if
and only if $\coker\alpha$ is torsion. If $s\in\mathcal{C}_R,$
$R/sR$ is torsion. Therefore the map $\alpha_s\colon R\rightarrow R$,
defined by $r\mapsto sr,$  is a full morphism and
$\alpha_s\otimes 1_B$ is invertible. Hence $b(s)$ is invertible in
$B.$ By the universal property of Ore localization there exists a
unique morphism of rings $\psi\colon A\rightarrow B$ such that
$\psi_{\mid_R}=b.$ Thus condition (ii) in
Theorem~\ref{def:universallocalization}  is satisfied.
\end{pro}

The following result is a generalization of
\cite[Theorem~1.3]{EisenbudRobson} to the semisimple situation. We
prove it using the theory of rank functions.

\begin{theo} \label{theo:essentialfinitelength}
Let $R$ be a hereditary noetherian semiprime ring which is a right
order in the semisimple ring $A$. Suppose there is no primitive central idempotent $e$ of $A$ such
that $eRe$ is simple artinian.  Let
$J\subseteq I$ be right ideals of $R$. Then $I/J$ is an artinian
right $R$-module if and only if $J$ is an essential submodule of
$I.$
\end{theo}
\begin{pro}
 We use the following known fact, see for example
\cite[Proposition~2.2.2]{Jategaonkarbook}: If $R$ is a right order
in a semisimple  ring, then a submodule $N$ of a torsion-free
right $R$-module $M$  is essential in $M$ if and only if $M/N$ is
torsion.

Suppose $I/J$ is artinian. Then it has finite length, that is, it
is a finite extension of simple right $R$-modules, and hence
torsion right $R$-modules by
Lemma~\ref{lem:semisimplerightordersandsimplemodules}. Thus $I/J$
is torsion.

On the other hand suppose $J$ is an essential submodule of $I.$ By
the remark at the beginning of the proof,  $I/J$ is a finitely
presented torsion right $R$-module. By the discussion on
Definition~\ref{de:exact} and
Proposition~\ref{theo:semihereditarytilting} it follows that $I/J$
has finite length.
\end{pro}

Now we come to the main result of this section.

\begin{theo}\label{theo:tiltingprimenoetherianrings}
Let $R$ be a hereditary noetherian prime ring which is not simple
artinian. Let $A_R$ be the simple artinian quotient ring of $R$.
Then
\begin{enumerate}[(1)]
\item $T=A\oplus A/R$ is a tilting right   $R$-module
with
$T^\perp=\mathcal{U}_r^\perp=\mathcal{V}_r^\perp=\mathcal{D}_r^\perp$.
\item $T=A\oplus A/R$ is a tilting left $R$-module
with
$T^\perp=\mathcal{U}_l^\perp=\mathcal{V}_l^\perp=\mathcal{D}_l^\perp.$
\item For any overring $R<S<A$ there exists a unique subset
$\mathcal{U}_S$ of $\mathcal{U}_r$ (respectively, of $\mathcal{U}_l)$ such that
$S\oplus S/R$ is a tilting right (left) $R$-module with tilting
class $\mathcal{U}_S^\perp.$ \item For any right Ore subset
$\mathfrak{S}$ of $\mathcal{C}_R,$ let
$$\mathcal{U}_{\mathfrak{S}}=\{U\in\mathcal{U}_r\mid
\textrm{for each } v\in U \textrm{ there exists } s\in\mathfrak{S} \textrm{
with } vs=0\}.$$ Then $R\mathfrak{S}^{-1}$ is the universal
localization of $R$ at $\mathcal{U}_\mathfrak{S}$. Moreover,
$T_{\mathfrak{S}}= R\mathfrak{S}^{-1}\oplus R\mathfrak{S}^{-1}/R$
is a tilting right $R$-module with tilting class
$T_{\mathfrak{S}}^\perp=\mathcal{U}_\mathfrak{S}^\perp$,
 and  $T_{\mathfrak{S}}$ is a
tilting left $R$-module with tilting class
$T_{\mathfrak{S}}^\perp=\{R/Rs\mid s\in\mathfrak{S}\}^\perp=\{\Tr
U\mid U\in \mathcal{U}_{\mathfrak{S}}\}^\perp$. \item For any (two
sided) Ore subset $\mathfrak{S}$ of $\mathcal{C}_R,$ let
$\mathcal{U}_{\mathfrak{S}}$ be as in (4). Then $T_{\mathfrak{S}}=
R\mathfrak{S}^{-1}\oplus R\mathfrak{S}^{-1}/R$ is a tilting right
$R$-module with tilting class $T_{\mathfrak{S}}^\perp=\{R/sR\mid
s\in\mathfrak{S}\}^\perp= \mathcal{U}_{\mathfrak{S}}^\perp.$
\end{enumerate}
\end{theo}

\begin{pro}
(1) By Proposition~\ref{theo:semihereditarytilting} and
Theorem~\ref{theo:universallocalizationwithfiltration},
$A/R$ is a directed union of
modules $N_i$
  where
each $N_i$ is a finite extension of simple right
$R$-modules. Hence, by
Corollary~\ref{theo:tiltingwithclass}, we get that $T$ is a
tilting module with tilting class
$T^\perp=\mathcal{U}_r^\perp.$ Moreover,
$\mathcal{U}_r^\perp=\mathcal{V}_r^\perp$ since  every
element in $\mathcal{V}_r$ is a finite extension of
elements in $\mathcal{U}_r$, see the discussion on
Definition~\ref{de:exact}. On the other hand, by
Corollary~\ref{coro:Oretiltingclass}, the tilting class of
$\mathcal{C}^{-1}R\oplus \mathcal{C}^{-1}R/R=T$ is
$\mathcal{D}_r^\perp.$

(2) is proven with symmetric arguments.

(3) It is proved in  \cite[Remark~3.3]{CrawleyBoeveytubes} that
every ring $S$ with $R<S<A$ is the universal localization of $R$
at some morphisms between finitely generated projective right
(left) $R$-modules. Recall from \cite[Chapter 5]{Schofieldbook}
that the universal localizations of $R$ embedding in $R_u$ are in
bijective correspondence with collections of stable association
classes of atomic full morphisms. So,  $S$ is the universal
localization of $R$ at a unique subset $\mathcal{U}_S$ of
$\mathcal{U}_r$ $(\mathcal{U}_l)$. Now, because of Proposition
\ref{theo:semihereditarytilting} (2) and (3),
  we can  apply
Corollary~\ref{coro:rhotorsiongivetilting}(2).

(4) $T_{\mathfrak{S}}$ is a  tilting left $R$-module with
$T_{\mathfrak{S}}^\perp=\{R/Rs\mid s\in\mathfrak{S}\}^\perp$ by
(the right version of) Corollary~\ref{coro:Oretiltingclass}.
Suppose  we have proved that
 $R\mathfrak{S}^{-1}$ is the universal
localization of $R$ at $\mathcal{U}_{\mathfrak{S}}.$  By
Proposition~\ref{theo:semihereditarytilting}, we can apply
Corollary~\ref{coro:rhotorsiongivetilting}(2)
 to obtain the desired results.

We now prove that $R\mathfrak{S}^{-1}$ is the universal
localization of $R$ at $\mathcal{U}_{\mathfrak{S}}$. The argument is very
similar to the one of \cite[Lemma~3.4]{CrawleyBoeveytubes}.

First of all, notice that ${}_RR\mathfrak{S}^{-1}$ is flat, and
for every $U\in \mathcal{U}_{\mathfrak{S}},$ $U\otimes_R
R\mathfrak{S}^{-1}=0.$ Hence, if $\exseqmap PQU\alpha$ is a
projective presentation of $U$ with $P$ and $Q$ finitely
generated, then $\alpha\otimes 1_{R\mathfrak{S}^{-1}}$ is
invertible. Hence condition (i) in
Theorem~\ref{def:universallocalization} is satisfied.

Let $B$ be a ring with a map $b\colon R\rightarrow B$ such that
for every $U\in\mathcal{U}_{\mathfrak{S}}$ and any finite
projective presentation $\exseqmap PQU\alpha,$
$\alpha\otimes_R1_B$ becomes invertible. Let $s\in\mathfrak{S}.$
Consider $R/sR.$ Since $\mathfrak{S}$ is right Ore, by
Theorem~\ref{theo:essentialfinitelength},  $R/sR$ has finite
length and therefore it has a finite filtration of simple right
$R$-modules. Recall $\mathfrak{S}$ consists of non-zero-divisors.
Hence $R/sR\cong s^{-1}R/R.$ Since $\mathfrak{S}$ is a right Ore
set, for every $s^{-1}r\in s^{-1}R,$ there exist
$t\in\mathfrak{S},\ x\in R$ such that $s^{-1}r=xt^{-1}.$
Therefore, for every $z\in R/sR$ there exists $t\in\mathfrak{S}$
with $zt=0.$ This implies that all the composition factors of
$R/sR$ are in $\mathcal{U}_\mathfrak{S}.$

For each $s\in\mathfrak{S},$ define the morphism $\delta_{s}\colon
R\rightarrow R,$ given by $r\rightarrow sr.$ By the foregoing,
 $\delta_s\otimes 1_B$ is invertible for every $s\in\mathfrak{S}$.
 Notice $\delta_s\otimes 1_B$ can be regarded as the morphism $B\rightarrow
 B$ defined by $x\mapsto b(s)x.$
Thus  $b(s)$ is invertible for all $s\in\mathfrak{S}.$ By the
universal property of Ore localization there exists a morphism of
rings $\gamma\colon R\mathfrak{S}^{-1}\rightarrow B$ making the
following diagram commutative
$$\xymatrix@C=0.3cm@R=0.3cm{ R\ar[rr]\ar[dr]_b & & R\mathfrak{S}^{-1}\ar[dl]^\gamma\\
& B }$$
Therefore condition (ii) in Theorem~\ref{def:universallocalization} is
satisfied.

For (5) apply the left and the right versions of (4).
\end{pro}

Before stating the next results we recall the following
definitions.

\begin{des}
(a) A right $R$-module $M$ is \emph{faithful} if the ideal $\ann
(M)=\{r\in R\mid mr=0\textrm{ for all } m\in M\}=0.$ We say $M$ is
\emph{unfaithful} if $M$ is not a faithful module.

(b) Let $Z$ be a commutative noetherian domain with quotient field
$K$ and let $Q$ be a central simple $K$-algebra. A
\emph{$Z$-order} in $Q$ is a $Z$-subalgebra $R$ of $Q$, finitely
generated as $Z$-module and such that $R$ contains a $K$-basis of
$Q$. A \emph{hereditary order} $R$ is a hereditary ring $R$ which
is a $Z$-order in some central simple $K$-algebra $Q$, where $Z$
is some Dedekind domain with quotient field $K\neq Z.$ The
hereditary order $R$ is a \emph{maximal order} if it is not
properly contained in any other $Z$-order in $Q$.
\end{des}

\begin{theo}\label{theo:alltiltingmaximalorders}
Let $R$ be a hereditary noetherian prime ring which is not
simple artinian. Let $\mathcal{U}_r$ be a set of
representatives of all isomorphism classes of all simple
right $R$-modules. Suppose that there are no simple
faithful right $R$-modules and that $\Ext_R^1(U_1,U_2)=0$
for any two non-isomorphic simple right $R$-modules
$U_1,U_2.$ Then
$$\mathbb{T}=\{T_\mathcal{W}=R_\mathcal{W}\oplus R_\mathcal{W}/R\mid
\mathcal{W}\subseteq\mathcal{U}_r\}$$ is a representative
set up to equivalence of the class of all tilting right
$R$-modules.

In particular, the statement holds true when $R$ is a maximal
order or a hereditary local noetherian prime ring which is not a
simple artinian ring (for example a not necessarily commutative
discrete valuation domain).
\end{theo}

\begin{pro}
For the first part we follow the terminology of \cite{Levysurvey}.
 If $M$ is a
finitely generated right $R$-module, then $M=P\oplus V$ where $P$
is projective and $V$ is the submodule consisting of the torsion
elements \cite[Lemma~5.7.4]{McConnellRobson}. Moreover, $V$ is of
finite length and has a decomposition $V=V_1\oplus V_2$ where
$V_1$ is a (finite) direct sum of uniserial modules whose
composition factors are all unfaithful, and $V_2$ is a direct sum
of modules whose composition factors belong to so-called faithful
towers \cite[Theorem~2.19]{Kuzmanovichhnprings} or
\cite[Theorem~4.6]{KlinglerLevy}. Since we are assuming that there
are no faithful simple right $R$-modules, $V_2=0.$

So $M^\perp=(W_1\oplus\dotsb\oplus W_n)^\perp$ where $W_i$
are indecomposable finitely generated uniserial modules.
Since $\Ext_R^1(U,U')=0$ for any two non-isomorphic simple
modules $U,U'$, we obtain that all composition factors of
$W_i$ are isomorphic to the same simple module $U_i.$ Hence
$W_i^\perp=U_i^\perp.$ Therefore
$M^\perp=U_1^\perp\cap\dotsb\cap U_n^\perp$. So for every
set of finitely generated right $R$-modules $\mathcal{V},$
there exists a set of simple right $R$-modules
$\mathcal{W}$ such that
$\mathcal{V}^\perp=\mathcal{W}^\perp.$ From
Remark~\ref{rems:boundmodules},
Proposition~\ref{theo:semihereditarytilting} and
Corollary~\ref{coro:rhotorsiongivetilting}(2) we infer that
$\mathbb{T}$ is a representative set  up to equivalence of
the class of all tilting right $R$-modules.

Let $R$ be a maximal order. That there are no faithful simple
$R$-modules is a concatenation of the following results. By
Lemma~\ref{lem:simplearetorsion}, every maximal right ideal
contains a nonzero divisor, so it is essential
\cite[2.1.15]{Jategaonkarbook}. Then the result follows by
\cite[Propositions~8.1, 8.3]{GoodearlWarfieldbook}. That maximal
orders verify the hypothesis on the extensions of simple modules
follows from \cite[Theorem~11.20]{GoodearlWarfieldbook} and
\cite[Theorem~22.4]{Reinerbook}.

If $R$ is a  hereditary local noetherian prime ring which is not a
simple artinian ring there is only one simple right $R$-module up
to isomorphism, and it is unfaithful since it is isomorphic to the
quotient of $R$ by its maximal ideal.
\end{pro}

We now recover the classification of tilting modules over Dedekind domains obtained in \cite[Theorem~5.3]{BazzoniEklofTrlifaj}.

\begin{coro}\label{coro:alltiltingDedekind}
Let $R$ be a Dedekind domain. \begin{enumerate}[(1)] \item  Let
$\mathfrak{M}$ be a subset of $\maxspec(R)$. Consider the
multiplicative subset of $R$
$\mathfrak{S}=R\setminus\operatornamewithlimits{\cup}\limits_{\mathfrak{p}\in\mathfrak{M}}\mathfrak{p}$
and the set of simple $R$-modules
$\mathcal{U}_{\mathfrak{S}}=\{R/\mathfrak{m}\mid
\mathfrak{m}\nsubseteq
\operatornamewithlimits{\cup}\limits_{\mathfrak{p}\in\mathfrak{M}}\mathfrak{p}\}$
(if $\mathfrak{M}=\emptyset,$ then $\mathfrak{S}=R\setminus
\{0\}$). Then $R\mathfrak{S}^{-1}$ is the universal localization
of $R$ at $\mathfrak{U}_{\mathfrak{S}}$ and
$T_{\mathfrak{G}}=R\mathfrak{G}^{-1}\oplus R\mathfrak{G}^{-1}/R$
is a tilting $R$-module with
$T_{\mathfrak{G}}^\perp=\mathcal{U}_{\mathfrak{S}}^\perp=\{R/sR\mid
s\in\mathfrak{G}\}^\perp$.

\item  Let
$\mathfrak{P}$ be a subset of $\maxspec(R)$. Consider
$\mathfrak{U}_\mathfrak{P}=\{R/\mathfrak{m}\mid
\mathfrak{m}\in\mathfrak{P}\}.$ Then
$T_{\mathfrak{P}}=R_{\mathcal{U}_{\mathfrak{P}}}\oplus
R_{\mathcal{U}_{\mathfrak{P}}}/R$ is a tilting right $R$-module
with $T_{\mathfrak{P}}^\perp=\mathcal{U}_{\mathfrak{P}}^\perp.$
Therefore the set $\{T_{\mathfrak{P}}\mid \mathfrak{P}\subseteq
\maxspec(R) \}$ is a representative set up to equivalence of the
class of all tilting $R$-modules.
\end{enumerate}
\end{coro}

\begin{pro}
(1) Recall that a Dedekind domain is a commutative noetherian prime
ring which is not simple artinian. Now notice that
$$\mathcal{U}_{\mathfrak{S}}=\{R/\mathfrak{m}\mid
\mathfrak{m}\nsubseteq
\operatornamewithlimits{\cup}\limits_{\mathfrak{p}\in\mathfrak{M}}\mathfrak{p}\}=
\{R/\mathfrak{m}\mid \textrm{ for every } v\in R/\mathfrak{m}
\textrm{ there is } s\in \mathfrak{S} \textrm{ with } vs=0\}.$$ Then apply
Theorem~\ref{theo:tiltingprimenoetherianrings}(5).

(2) By Proposition~\ref{theo:semihereditarytilting}
and Corollary~\ref{coro:rhotorsiongivetilting}(2), we obtain that
$T_{\mathfrak{P}}$ is a tilting module with
$T_{\mathfrak{P}}^\perp=\mathcal{U}_{\mathfrak{P}}^\perp.$

For the second statement, we prove that $R$ satisfies the
conditions of Theorem~\ref{theo:alltiltingmaximalorders}.
All simple $R$-modules are unfaithful because they are
isomorphic to the quotient of $R$ by a maximal ideal. On
the other hand let $M$ be an extension of the nonisomorphic
simple $R$-modules $R/\mathfrak{p}$ and $R/\mathfrak{q}$
with $\mathfrak{p}$ and $\mathfrak{q}$ nonzero prime ideals
of $R$. Notice that $\ann(M)=\mathfrak{p}\mathfrak{q}.$
Therefore $M\cong T_{\mathfrak{p}}(M)\oplus
T_\mathfrak{q}(M)$ where $T_{\mathfrak{p}}(M)$ and
$T_\mathfrak{q}(M)$  denote the $\mathfrak{p}$-primary and
the $\mathfrak{q}$-primary components of $M.$ Therefore
$M\cong R/\mathfrak{p}\oplus R/\mathfrak{q}$, see
\cite[Sections~5.1, 6.3]{BerrickKeating}.
\end{pro}

\medskip


{\bf Acknowledgments}\medskip

We would like to thank Dolors Herbera for interesting
conversations and valuable suggestions. Among other things,
she pointed out to us Remarks~\ref{rems:boundmodules}
and~\ref{rems:universalisepi}(4), and she suggested the proof of
Theorem~\ref{theo:alltiltingmaximalorders} for Dedekind domains.

This research   was carried out during two visits of the second named author at Universit\`a dell'Insubria, Varese,
 in 2005 and 2006, supported by  a  grant of the Facolt\`a di Scienze dell'Universit\`a dell'Insubria, Varese,
 and by Departament d'Universitats, Recerca i Societat de la Informaci\'o de la Generalitat de
 Catalunya. The second named author would like to thank Dipartimento di Informatica
e Comunicazione  dell'Universit\`a degli Studi dell'Insubria di
Varese for its hospitality.

   First named author partially supported by PRIN 2005 ``Prospettive in
teoria degli anelli, algebre di Hopf e categorie di moduli''.\\
Both authors partially supported by the DGI and the European
Regional Development Fund, jointly, through Project
 MTM2005--00934, and by the Comissionat per Universitats i Recerca
of the Generalitat de Ca\-ta\-lunya, Project 2005SGR00206.


\bibliographystyle{amsalpha}
 \bibliography{grupitosbuenos}

 \end{document}